%% file: GAPManifolds.tex
\documentclass[english]{jnsao}

\usepackage[margin=1.3in]{geometry}

\usepackage{epsfig}
\usepackage{booktabs}
\usepackage{enumerate}
\usepackage{url}
\usepackage{euscript}

\newtheorem{thm}{Theorem}[section]
\newtheorem{fact}{Fact}[section]
\newtheorem{prp}{Proposition}[section]

\newtheorem{ass}{Assumption}[section]
\newtheorem{lem}{Lemma}[section]
\newtheorem{rem}{Remark}[section]

\newtheorem{defin}{Definition}[section]
\newtheorem{ex}{Example}[section]
\newenvironment{pf}{\smallbreak\noindent{\textit{Proof. }}}{\hfill$\Box$\smallbreak}

\makeatletter\@addtoreset{case}{thm}\makeatother
\makeatletter\@addtoreset{case}{subsection}\makeatother

\makeatletter\@addtoreset{subcase}{case}\makeatother

\manuscriptsubmitted{2021-01-28}
\manuscriptaccepted{2023-02-21}
\manuscriptvolume{4}
\manuscriptnumber{7139}
\manuscriptyear{2023}
\manuscriptdoi{10.46298/jnsao-2023-7139}

  \manuscriptlicense{CC-BY 4.0}

  \shortauthor{Fält, Giselsson}

 \shorttitle{GAP on Manifolds and Convex Sets}

\usepackage{tikz}
\usetikzlibrary{calc}
\usetikzlibrary{patterns}
\usetikzlibrary{arrows.meta}
\usepackage[utf8]{inputenc}
\usepackage{pgfplots}
\pgfplotsset{compat=1.17}

\usepackage{mathtools}
\usepackage{multicol}
\usepackage{cases}
\usepackage{enumitem}

\newcommand{\argmin}{\mathop{\mathrm argmin}}

\newcommand{\fix}{{\mathrm{fix}}}

\newcommand{\uset}{{\mathcal{U}}}
\newcommand{\vset}{{\mathcal{V}}}
\newcommand{\reals}{\mathbb{R}}
\newcommand{\ldef}{\coloneqq}
\newcommand{\rdef}{\eqqcolon}
\newcommand{\re}{\text{\normalfont Re}\,}

\newcommand{\sign}[1]{\text{\normalfont sgn}(#1)}
\newcommand{\bd}{\text{\normalfont bd}\,}
\newcommand{\inter}{\text{\normalfont int}}

\newcommand{\realR}{\mathbb{R}} 
\newcommand{\naturalN}{\mathbb{N}}

\newcommand{\M}{{\mathcal{M}}}
\newcommand{\N}{{\mathcal{N}}}
\newcommand{\cset}{{\mathcal{C}}}
\renewcommand{\C}{{\mathcal{C}}}
\newcommand{\Proj}{\Pi}
\newcommand{\Normal}{{\textrm{N}}}
\newcommand{\T}{\mathrm{T}}
\newcommand{\Ball}{{\mathcal{B}}}
\newcommand{\Ballo}{\mathcal{B}^o}

\newcommand{\mnote}[1]{}

\newcommand{\fnote}[1]{}

\newcommand{\mnoteh}[1]{}
\newcommand{\cmt}[1]{}

\newcommand{\norm}[1]{\left\lVert#1\right\rVert}

\newcommand{\spectrum}[1]{\Lambda(#1)}
\newcommand{\jacobian}[2]{\mathrm{J}_{#1}(#2)}
\newcommand{\sreg}{\mathrm{sr}}
\newcommand{\reg}{\mathrm{r}}
\newcommand{\dist}{\textrm{d}}
\newcommand{\rate}[1]{\sigma(#1)}

\newcommand{\sinS}{\EuScript{S}}
\newcommand{\cosC}{\EuScript{C}}

\begin{document}
\title{Generalized Alternating Projections on Manifolds and Convex Sets}
\author{Mattias Fält%
\thanks{Department of Automatic Control, Lund University, Sweden. \email{mattias.falt@control.lth.se}.}, Pontus Giselsson\thanks{Department of Automatic Control, Lund University, Sweden. \email{pontus.giselsson@control.lth.se}.}}
\maketitle

\begin{abstract}
In this paper, we extend the previous convergence results for the generalized alternating projection method applied to subspaces in \cite{falt-optimal} to hold also for smooth manifolds.
We show that the algorithm locally behaves similarly in the subspace and manifold settings and that the same rates are obtained.
We also present convergence rate results for when the algorithm is applied to non-empty, closed, and convex sets.
The results are based on a finite identification property that implies that the algorithm after an initial identification phase solves a smooth manifold feasibility problem. Therefore, the rates in this paper hold asymptotically for problems in which this identification property is satisfied. We present a few examples where this is the case and also a counter example for when this is not.
\end{abstract}

\mnoteh{Papers to look at: About Transversality of collections of sets.}
\section{Introduction}
The problem of finding a point in the intersection of sets has a long history with many proposed algorithms.
They generally rely on successive projections onto the respective sets.
The method of alternating projections (MAP, or AP) was famously studied by von Neumann \cite{vonNeumann} for the case of two subspaces
and has a wide range of applications \cite{Deutsch1992}.
Many variants have been suggested and shown to converge in the case of convex sets,
for example using relaxed projections \cite{GAP_Agmon, GAP_Motzkin, GAP_Bregman,GAP_Gubin} or inexact projections \cite{kruger2016regularity},
Dykstra's algorithm~\cite{Boyle1986}, Douglas--Rachford splitting~\cite{DouglasRachford,LionsMercier1979},
and its dual algorithm ADMM~\cite{Glowinski1975,BoydDistributed}.

Many results on the linear convergence rates of these algorithms have been shown
and are generally stated as a function of a regularity constant such as the smallest angle between the sets,
which in the case of affine sets is known as the Friedrichs angle $\theta_F$.
In the case of two subspaces, the method of alternating projections was shown to converge with
the linear rate $\cos^2(\theta_F)$ \cite{Deutsch1995},
and the Douglas--Rachford method with the rate $\cos(\theta_F)$ \cite{Bauschke_lin_rate_Friedrich}.
In \cite{Bauschke_opt_rate_matr}, the authors studied a few methods with relaxed projections
and the optimal rates with respect to the relaxation parameters were found.
The generalized alternating projection (GAP)---which generalizes most of the algorithms above by allowing several relaxation parameters---was studied in \cite{falt-optimal,dao2018linear,dao2019linear}.
It was shown in \cite{falt-optimal} that the faster rate $\frac{1-\sin(\theta_F)}{1+\sin(\theta_F)}$ is achievable with the right parameters.
It was also shown that, under general assumptions,
this is the best possible rate for this generalization.

When it comes to general convex sets,
local linear convergence of these algorithms is not guaranteed.
Several different assumptions on the intersection between the sets have been proposed and shown to be sufficient.
Some of these assumptions include linear regularity or bounded linear regularity, see for example \cite{lewis2009local, bauschke1993convergence}.
An overview of set regularities can be found in \cite{kruger2006regularity} and a survey on their relation can be found in \cite{kruger2018set}.
\mnoteh{MORE RESULTS ON CONVEX}
Under subtransversality assumptions of two convex sets, the R-linear rate presented in \cite{Luke2020} translates to a $\cos(\theta_F/2)$ contraction rate for the Douglas--Rachford algorithm, when translated to the subspace setting.

For general non-convex sets, convergence to a feasible point cannot be guaranteed,
and local convergence is instead studied.
For the alternating projections method, different types of regularity have been shown to be sufficient
for local linear convergence \cite{lewis2009local,Bauschke13-theory,Bauschke13-applications,Noll13}.
\mnoteh{Noll :General sets, local convergence if separable intersection}
For the alternating projections algorithm, the results in \cite{lewis2009local} for possibly non-convex super-regular sets with linearly regular intersection translates to the known optimal rate of $\cos^2(\theta_F)$ when applied to sub-spaces.
In \cite{Drusvyatskiy15}, the authors showed that a transversality property can be used to guarantee local linear convergence.
However, both the assumptions and rates presented in this paper are quite conservative. For example, in the case of two subspaces, the rate presented in \cite{Drusvyatskiy15} translates to $\cos^2(\theta_F/2)$
which is considerably worse than the known contraction rate $\cos(\theta_F)$ and the local linear rate $\cos^2(\theta_F)$.
Among the few known results for alternating relaxed projections,
local linear convergence was shown for the MARP algorithm in \cite{bauschke2014MARP} under different regularity assumptions. However, this paper assumes that the projections are under-relaxed,
which was shown in \cite{falt-optimal} to result in sub-optimal local rates.

One approach to show local convergence rates for general convex sets is by showing that the algorithms
eventually project onto subsets that have nicer properties,
i.e., that the algorithm identifies these subsets in finite time.
This can be done by partitioning the boundary of sets into a collection of smooth manifolds,
and then studying the algorithm on these manifolds.
There has been a lot of research into these identification properties for various algorithms,
see for example \cite{hare-2004,Lewis2011Activity,act_ident_15}.
However, as far as the authors know, none of these results apply to projection methods on feasibility problems.
The fundamental problem seems to be that gradients are vanishing at any feasible point when a feasibility problem is reformulated as an optimization problem, so the regularity assumptions are therefore not satisfied.
However, for specific problems it can sometimes be known that the algorithm will identify such surfaces,
for example when the entire boundary is a smooth manifold,
or when the algorithm is known to converge to the relative interior of one of the manifolds.

\mnote{Something about problems posed with manifolds}

In this paper, we study generalized alternating projections in the setting of two smooth manifolds. The special case of alternating projections is studied in \cite{lewis-2008}. There, it is shown that the smooth manifolds locally can be approximated by affine sets and that the convergence rates known from affine sets translate to local linear rates in this setting under a transversality condition.
A similar result is found in \cite{andersson2013alternating} under slightly relaxed assumptions. We show that the weaker assumption in \cite{andersson2013alternating}
is sufficient to show local linear convergence also of the generalized alternating projections method on smooth manifolds.
Moreover, we show that the optimal rates and parameters for linear subspaces found in \cite{falt-optimal}
translate to the smooth manifold setting.

We combine our rate results for generalized alternating projections on smooth manifolds with a finite identification property. This gives convergence rate results for the algorithm when applied to convex sets for which the algorithm enjoys this identification property. We provide some classes of convex sets for which this property holds, implying that the convergence rate result for manifolds and subspaces is valid also for these sets. We also provide one counter-example where we illustrate
that even in the setting of polyhedral sets and the presence of regularity,
the finite identification property does not hold. As a consequence, the problem can in that case not be locally reduced to that of affine sets, as is the case for alternating projections. \mnote{CITE, Not really what I want to say, but (linear convergence using linear regularity through, e.g., \cite{lewis2009local}, and \cite{bauschke1996projection} 5.7.2 to show linear Regularity)}

\mnote{
%
%
%
%
}

\section{Notation}
We let $\naturalN$ denote the set of non-negative integers, $\reals$ be the real line, $\reals^n$ be the set of $n$-dimensional real vectors, and $\reals^{n\times m}$ be the set of $n\times m$ real matrices.  
We denote the identity operator by $I$ and the operator norm by $\|\cdot\|$.
For a matrix $A\in\mathbb{R}^{n\times n}$ we let $\spectrum{A}$ be the set of eigenvalues and $\rho(A)\ldef \max_{\lambda\in\spectrum{A}}|\lambda|$ the spectral radius.
If the limit $\lim_{k\rightarrow\infty} A^k$ exists, we denote it by $A^\infty$ and define $\rate{A}\ldef\|A-A^\infty\|$. For a vector $v\in\realR^n$ we also denote the vector norm by $\|v\|\ldef\sqrt{\langle v, v\rangle}$.
The Jacobian of a function $F:\reals^m\to\reals^n$ at a point $x\in\reals^n$ is denoted by $\jacobian{F}{x}$.
We denote the closed ball around a point
$x\in\realR^n$ and with radius $\delta> 0$, i.e., $\{y\in\realR^n\mid \|x-y\|\leq\delta\}$, by $\Ball_\delta(x)$
and the open ball $\{y\in\realR\mid \|x-y\|<\delta\}$ by $\Ballo_\delta(x)$.
\mnote{Differentiate between Identity operator and matrix.}
\mnote{$\realR$ or euclidean space?}

\section{Preliminaries}
\begin{defin}[projection]
    The projection of an element $x\in\realR^n$ onto a closed, non-empty subset $C\subset \realR^n$ is defined by
    \begin{align*}
        \Proj_C(x) \ldef  \argmin_{y\in C} \|x-y\|
    \end{align*}
    when the argmin is unique.
\end{defin}
\begin{defin}[relaxed projection]
  The relaxed projection of an element $x\in\realR^n$ onto a closed, non-empty subset $C\subset \realR^n$ with relaxation parameter $\alpha\neq 0$ is defined as
    \begin{align*}
        \Proj_{C}^{\alpha}(x)\ldef (1-\alpha)x+\alpha \Proj_{\cset}(x)
    \end{align*}
    when the argmin is unique.
\end{defin}
\subsection{Subspaces}
In this section, we introduce some basic properties of subspaces that will be useful in the study of the local properties of manifolds.

\begin{defin}\label{def:principal}
The \emph{principal angles} $\theta_{k}\in[0,\pi/2],\,k=1,\dots,p$ between two subspaces
$\uset,\vset\subset\realR^n$, where $p=\min(\dim\uset,\dim\vset)$, are recursively defined by
\begin{eqnarray*}
\cos\theta_{k} & \ldef  & \max_{u_{k}\in\uset,\,v_{k}\in\vset}\left\langle u_{k},v_{k}\right\rangle \\
 & \text{s.t.} & \norm{u_{k}}=\norm{v_{k}}=1,\\
 &  & \left\langle u_{k},v_{i}\right\rangle =\left\langle u_{i},v_{k}\right\rangle =0,\forall\,i=1,\ldots,k-1.
\end{eqnarray*}
\end{defin}
\begin{fact}\label{fct:friedrichs}
\cite[Def. 3.1, Prop. 3.3]{Bauschke_opt_rate_matr}
The principal angles
are unique and satisfy $0\leq\theta_{1}\leq\theta_{2}\leq\dots\theta_{p}\leq\pi/2$.
The angle $\theta_{F}\ldef \theta_{s+1}$, where $s=\text{dim}(\uset\cap\vset)$, is the \emph{Friedrichs angle} and it
is the smallest non-zero principal angle.
\end{fact}

The cosine of the Friedrichs angle occurs naturally in many convergence rate results and is denoted as in the following definition.
\begin{defin}\label{defin:friedrichs}
  The cosine of the Friedrichs angle $\theta_F$ between two subspaces $\,\uset,\vset\subset\realR^n$ is denoted as
    \begin{align*}
        c(\uset, \vset) := \cos(\theta_F).
    \end{align*}
\end{defin}

We see that $\theta_i=0$ if and only if $i\leq s$, where $s=\text{dim}(\uset\cap\vset)$, so $\theta_F$ is well defined whenever $\min(\dim\uset,\dim\vset)=p>s=\dim(\uset\cap\vset)$,
i.e., when no subspace is contained in the other.
\mnoteh{When one subspace is contained in the other we define $c(\uset,\vset)=0$.}

\begin{defin}\label{def:linear}
$A\in\realR^{n\times n}$ is $\emph{\text{linearly convergent}}$
to $A^{\infty}\in\realR^{n\times n}$ with \emph{linear convergence rate }$\mu\in[0,1)$
if there exist $M,N>0$ such that
\[
\norm{A^{k}-A^{\infty}}\leq M\mu^{k}\quad\forall k>N,\,k\in\naturalN.
\]
\end{defin}
\begin{defin}\label{def:limit}
\cite[Fact 2.3]{Bauschke_opt_rate_matr}  For $A\in\realR^{n\times n}$ we say that
$\lambda\in\spectrum{A}$ is \emph{semisimple} if $\text{ker}(A-\lambda I)=\text{ker}(A-\lambda I)^{2}.$
\end{defin}
\begin{fact}\label{fct:limitexists}
\cite[Fact 2.4]{Bauschke_opt_rate_matr} For $A\in\realR^{n\times n}$, the limit
$A^{\infty}\ldef \lim_{k\rightarrow\infty}A^{k}$ exists if and only if
\end{fact}
\begin{itemize}
\item $\rho(A)<1$ or
\item $\rho(A)=1$ and $\lambda=1$ is semisimple and the only eigenvalue
on the unit circle.
\end{itemize}
\begin{defin}\label{def:subdominant}
\cite[Def. 2.10]{Bauschke_opt_rate_matr} Let $A\in\realR^{n\times n}$ be a matrix with $\rho(A)\leq 1$ and define
\[
\gamma(A)\ldef\max\left\{|\lambda|\, \mid\, \lambda\in\{0\}\cup\spectrum{A}\setminus\{1\}\right\} .
\]
Then $\lambda\in\spectrum{A}$ is a \emph{subdominant eigenvalue} if  $|\lambda|=\gamma(A)$.
\end{defin}%

\begin{fact}
\label{fact:Convergent-at-rate}\cite[Thm. 2.12]{Bauschke_opt_rate_matr}
If $A\in\realR^{n\times n}$ is convergent to $A^{\infty}\in\realR^{n\times n}$
then
\end{fact}
\begin{itemize}
\item $A$ is linearly convergent with any rate $\mu\in(\gamma(A),1)$
\item If $A$ is linearly convergent with rate $\mu\in[0,1)$, then $\mu\in[\gamma(A),1)$.
\end{itemize}

\subsection{Manifolds}

The following definitions and results follow those in \cite{lewis-2008}.

\begin{defin}[smooth manifold]
    A set $\M\subset\realR^n$ is a $\C^k$-manifold around a point $x\in\M$ if there is an open set $U\subset\realR^n$ containing $x$ such that
    \begin{align*}
        \M\cap U = \{ x : F(x)=0\}
    \end{align*}
    where $F:U\rightarrow \realR^d$ is a $\C^k$ function with surjective derivative throughout $U$.
\end{defin}

\begin{defin}[tangent space and tangent plane]
    The tangent space to a manifold $\M$ at $x\in\reals^n$ is given by
    \begin{align*}
        \T_\M(x) = \ker \jacobian{F}{x}.
    \end{align*}
and is independent of the choice of $F$ that defines the manifold. The {\emph{tangent plane}} is $\T_\M(x)+\{x\}$.
\end{defin}

\begin{defin}[normal vector]
    $v\in\realR^n$ is a normal vector to the manifold $\M$ at $x\in\realR^n$ if $\langle v, t\rangle=0$ for all $t\in\T_\M(x)$.
\end{defin}

\begin{defin}[smooth boundary]
    We say that a closed set $C\subset\realR^n$ has a $\C^k$ smooth boundary around $\bar{x}\in\realR^n$ if $\bd(C)$ is a $\C^k$ smooth manifold around $\bar{x}$.
\end{defin}

\begin{rem}\label{rem:smooth-normal}
We note that if a set $C\in\realR^n$ is solid, i.e., $\inter(C)\neq\emptyset$, with a $\C^k$ smooth boundary around some point $\bar{x}$, then the boundary is defined in some neighborhood $U$ of $\bar{x}$ by some $f:\realR^n\rightarrow\realR$ as $\bd(C)\cap U = \{x: f(x)=0\}$.
The tangent space given by $\ker{\jacobian{f}{x}}$ is therefore an $\realR^{n-1}$ dimensional plane, with normal vector $\nabla f(x)$. Since $f$ is a $\C^k$ smooth function, the normal vector is a $\C^{k-1}$ smooth function of $x$.
\end{rem}

We now define the regularity condition that will be sufficient to show linear convergence of the GAP method.

\begin{ass}[regularity]\label{ass:regularity}%
    Two manifolds $\M$ and $\N$ satisfy the regularity assumption at a point $x\in\reals^n$ if they are $\C^k$-smooth ($k\geq 2$) around $x\in\M\cap\N$ and
    \begin{enumerate}[label=A\arabic*.,ref=A\arabic*]
        \item $\M\cap\N$ is a $\C^k$ smooth manifold around $x$\label{ass:reg-smooth}
        \item $\T_{\M\cap\N}(x)=\T_\M(x)\cap \T_\N(x)$. \label{ass:reg-intersection}
    \end{enumerate}
\end{ass}

We note that for closed convex sets, the assumption \ref{ass:reg-intersection} is the conical hull intersection property \cite{chui1990constrained,bakan2005strong}. We also note that our regularity condition is equivalent to the one used in \cite{andersson2013alternating} for proving a linear convergence rate for alternating projections. Besides \ref{ass:reg-smooth}, they use the assumption that the manifolds are {\emph{non-tangential}}, as defined in \cite[Def.~3.4]{andersson2013alternating}, in place of \ref{ass:reg-intersection}. These two latter properties are equivalent except when one manifold locally is a subset of the other, which is an uninteresting trivial case for these methods, see \cite[Prop. 3.2 and Prop. 3.5]{andersson2013alternating}. Another common regularity property that has been used, e.g., in \cite{lewis-2008} is transversality.

\begin{defin}[transversality]\label{def:substransversal-manifolds} Two  $\mathcal{C}^k$-smooth manifolds $\M$ and $\N$ are transversal at $\bar{x}\in\realR^n$ if $\T_{\M}(\bar{x})+T_{\N}(\bar{x})=\realR^{n}$.
\end{defin}

We note that both~\ref{ass:reg-smooth} and~\ref{ass:reg-intersection} in Assumption~\ref{ass:regularity} are implied by the transversality assumption \cite{kruger2018set}.
However, transversality is not a consequence of Assumption \ref{ass:regularity} as we see in the following example.

\begin{ex}
Let $\M=\{(x,0,x^2)\mid x\in\realR\}$ and $\N=\{(0,y,0)\mid y\in\realR\}$ where $\M\cap\N=\{0\}$. We have $\T_\M(0)=\{(x,0,0)\mid x\in\realR\}$ and $\T_\N(0)=\N$.
So the manifolds clearly satisfy Assumption \ref{ass:regularity} at $0$, but not the transversality condition $\T_\M(0)+\T_\N(0)=\{(x,y,0)\mid x,y\in\realR\}\neq \realR^3$.
\end{ex}

With some abuse of notation, we define the angle between two manifolds at a point in their intersection
using their tangent spaces.
\begin{defin}\label{defin:cmn}
    For $x\in\M\cap\N$ let
    \begin{align*}
        c(\M,\N,x)\ldef c(\T_\M(x),\T_\N(x)).
    \end{align*}
\end{defin}

The regularity condition implies that both the manifolds and their intersection locally behave similarly to their tangent planes.
In particular, the angle between two lines that belong to different tangent planes is zero if and only if the lines are parallel to the intersection of the manifolds, as seen by \ref{ass:reg-intersection}.
This is crucial to show linear convergence.
We also note that, under the regularity assumptions, the Friedrichs angle $\theta_F$ is positive unless one manifold is locally a subset of the other.
To see this, we know that $\theta_F$ is well defined and positive unless one tangent space is a subset of the other,
for example $\T_\M(x)\subset \T_\N(x)$.
But since $\dim(\T_\M(x))=\dim(\M)$ around $x$, \ref{ass:reg-intersection} implies that also $\dim(\M)=\dim(\M\cap\N)$ around $x$,
i.e., that $\M$ locally is a subset of $\N$.
Under the regularity assumption, we therefore either have a positive Friedrichs angle
or a locally trivial problem.

Next, we show that relaxed projections are locally well defined on smooth manifolds,
and that their Jacobian is given by relaxed projections onto their tangent spaces.
By well defined we mean that the relaxed projection point exists and is unique.

The following Lemma is from \cite[Lem 4]{lewis-2008}.
\begin{lem}[projection onto manifold]\label{lem:proj-manifold}%
If $\M$ is a $\C^k$ manifold (with $k\geq 2$) around $\bar{x}\in \M$, then $\Proj_{\M}$ is well defined and $\C^{k-1}$ around $\bar{x}$. Moreover $\jacobian{\Proj_\M}{\bar{x}}=\Proj_{\T_{\M}(\bar{x})}$.
\end{lem}

\begin{lem}[relaxed projection onto manifold]\label{lem:proj-manifold-relaxed}%
If $\M$ is a $\C^k$ manifold (with $k\geq 2$) around $\bar{x}\in \M$, then $\jacobian{\Proj^{\alpha}_\M}{\bar{x}}=\Proj_{\T_{\M}(\bar{x})}^{\alpha}$,
    and $\Proj_{\M}^{\alpha}$
    are well defined and $\C^{k-1}$ around $\bar{x}$.

\end{lem}
\begin{pf}
    $\jacobian{\Proj_\M^{\alpha}}{\bar{x}}=\jacobian{(1-\alpha)I+\alpha \Proj_\M}{\bar{x}}=(1-\alpha)I+\alpha \Proj_{\T_{\M}(\bar{x})}=\Proj_{\T_{\M}(\bar{x})}^{\alpha}$.
    The result now follows from Lemma \ref{lem:proj-manifold}.
\end{pf}
%

\section{Generalized Alternating Projections}
In this section, we define the generalized alternating projections (GAP) operator, and state some known results. We denote the feasibility problem of finding $x\in\uset\cap\vset$ by $(\uset,\vset)$ to signify that the algorithm depends on the ordering of the two sets.

\begin{defin}[generalized alternating projections]\label{def:GAP}%
    The generalized alternating projections algorithm (GAP)~\cite{GAPLS} for the feasibility problem ($\uset,\vset$), where $\uset,\vset\subseteq\reals^n$ and $\uset\cap\vset\neq\emptyset$, is defined by the iteration
    \begin{equation}\label{eq:GAP}
    x_{k+1}\ldef Sx_{k},
    \end{equation}
    where
    \begin{equation}\label{eq:GAP2}
    S=(1-\alpha)I+\alpha \Proj_{\uset}^{\alpha_{2}}\Proj_{\vset}^{\alpha_{1}}=:\,(1-\alpha)I+\alpha T
  \end{equation}
  and $\alpha,\alpha_1,\alpha_2\in\reals$ are scalar parameters of the algorithm.
\end{defin}

For closed convex sets, the operator $S$ is averaged and the iterates converge to a point in the fixed-point set $\fix S$
under Assumption~\ref{ass:alpha}, see, e.g.,~\cite{GAPLS} where these results are collected.
\begin{ass}\label{ass:alpha}%
    Assume that $\alpha\in(0,1]$, $\alpha_{1},\alpha_{2}\in(0,2]$ and that one of the following holds
    \begin{enumerate}[label=B\arabic*.,ref=B\arabic*]
    \item $\alpha_{1},\alpha_{2}\in(0,2)$\label{ass:alpha-1}
    \item $\alpha\in(0,1)$ with either $\alpha_1\neq2$ or $\alpha_2\neq2$\label{ass:alpha-2}
    \item $\alpha\in(0,1)$ and $\alpha_{1}=\alpha_{2}=2$\label{ass:alpha-3}
    \end{enumerate}
\end{ass}
The following result is shown in \cite{GAPLS}.
\begin{lem}\label{lem:fixed-points}%
    Let $\uset,\vset\subset\reals^n$ be two non-empty linear subspaces and consider the feasibility problem $(\uset,\vset)$.
    The fixed point set $\fix S\ldef\{x\in\reals^n\mid Sx=x\}$ of the GAP operator $S$ in \eqref{eq:GAP} is: $\uset\cap\vset$ under Assumption \ref{ass:alpha} case \ref{ass:alpha-1} and \ref{ass:alpha-2},
    and $\,\uset\cap\vset + (\uset^\perp\cap\vset^\perp)$ under Assumption \ref{ass:alpha} case \ref{ass:alpha-3}.
\end{lem}

Throughout this section, we assume that the subspaces $\uset,\vset\subset\reals^n$ are non-empty, which implies that the problem $(\uset,\vset)$ is consistent and satisfies $0\in\uset\cap\vset$.

The following proposition and remark are found in~\cite[Prop. 3.4]{Bauschke_opt_rate_matr} and \cite{falt-optimal} respectively.
    \begin{prp}\label{prp:projections}%
      Let $\uset$ and $\vset$ be subspaces in $\realR^{n}$ and let $p\ldef \dim(\uset)$ and $q\ldef \dim(\vset)$ 
        satisfy
        $p\leq q$, $p+q< n$ and $p,q\geq1$.
        Then the projection matrices $\Proj_{\uset}\in\reals^{n\times n}$
        and $\Proj_{\vset}\in\reals^{n\times n}$ become

    \begin{align}
    \Proj_{\uset} & =D\begin{pmatrix}I_{p} & 0 & 0 & 0\\
    0 & 0_{p} & 0 & 0\\
    0 & 0 & 0_{q-p} & 0\\
    0 & 0 & 0 & 0_{n-p-q}
    \end{pmatrix}D^{*},\label{eq:Pu}\\
    \Proj_{\vset} & =D\begin{pmatrix}\cosC^{2} & \cosC\sinS & 0 & 0\\
    \cosC\sinS & \sinS^{2} & 0 & 0\\
    0 & 0 & I_{q-p} & 0\\
    0 & 0 & 0 & 0_{n-p-q}
    \end{pmatrix}D^{*}\label{eq:Pv}
    \end{align}
    and
    \begin{equation}
    \Proj_{\uset}\Proj_{\vset}=D\begin{pmatrix}\cosC^{2} & \cosC\sinS & 0 & 0\\
    0 & 0_{p} & 0 & 0\\
    0 & 0 & 0_{q-p} & 0\\
    0 & 0 & 0 & 0_{n-p-q}
    \end{pmatrix}D^{*},
    \end{equation}
    where $\cosC$ and $\sinS$ are diagonal matrices containing the cosine and sine
    of the principal angles $\theta_{i}$, i.e.,
    \begin{align*}
     \sinS&=\begin{bmatrix}\sin(\theta_{1})&&\\&\ddots &\\&& \sin(\theta_{p})\end{bmatrix},&
    \cosC&=\begin{bmatrix}\cos(\theta_{1})&&\\&\ddots&\\&&\cos(\theta_{p})\end{bmatrix},
    \end{align*}
    and $D\in\realR^{n\times n}$ is an orthogonal matrix.%
\end{prp}

Under the assumptions in Proposition~\ref{prp:projections},
the linear operator $T$ that is implicitly defined in~\eqref{eq:GAP2} becomes
\begin{align*}
T & =  \Proj_{\uset}^{\alpha_{2}}\Proj_{\vset}^{\alpha_{1}}=((1-\alpha_{2})I+\alpha_{2}\Proj_{\uset})((1-\alpha_{1})I+\alpha_{1}\Proj_{\vset})\\
 & =  (1-\alpha_{2})(1-\alpha_{1})I+\alpha_{2}(1-\alpha_{1})\Proj_{\uset}+\alpha_{1}(1-\alpha_{2})\Proj_{\vset}+\alpha_{1}\alpha_{2}\Proj_{\uset}\Proj_{\vset}\\
                  & =  D\begin{bmatrix}T_{1}&0&0\\0&T_{2}&0\\0&0&T_{3}\end{bmatrix}D^{*}
\end{align*}
where
\begin{align}
T_{1} & =\begin{pmatrix}I_{p}-\alpha_{1}\sinS^{2} & \alpha_{1}\cosC\sinS\\
\alpha_{1}(1-\alpha_{2})\cosC\sinS & (1-\alpha_{2})(I_{p}-\alpha_{1}\cosC^{2})
\end{pmatrix},\label{eq:T1}\\
T_{2} & =(1-\alpha_{2})I_{q-p},\quad T_{3}=(1-\alpha_{2})(1-\alpha_{1})I_{n-p-q}.\nonumber
\end{align}
The rows and columns of $T_1$ can be reordered so that it is a block-diagonal
matrix with blocks
\begin{equation}
T_{1_i}\hspace{-0.05cm}=\hspace{-0.05cm}\begin{pmatrix}1-\alpha_{1}s_{i}^{2} & \alpha_{1}c_{i}s_{i}\\
\alpha_{1}(1-\alpha_{2})c_{i}s_{i} & (1-\alpha_{2})(1-\alpha_{1}c_{i}^{2})
\end{pmatrix},\, i\in1,\dots,p\label{eq:T1matrix}
\end{equation}
where $s_{i}\ldef \sin(\theta_i),\,c_{i}\ldef \cos(\theta_{i})$.
The eigenvalues of $T$ are therefore
    $\lambda^{3}\ldef (1-\alpha_{2})$,
    $\lambda^{4}\ldef (1-\alpha_{2})(1-\alpha_{1})$,
and for every $T_{1_i}$
\begin{align}
  \lambda_{i}^{1,2} & =\frac{1}{2}\left(2-\alpha_{1}-\alpha_{2}+\alpha_{1}\alpha_{2}c_{i}^{2}\right)\pm\sqrt{\frac{1}{4}\left(2-\alpha_{1}-\alpha_{2}+\alpha_{1}\alpha_{2}c_{i}^{2}\right)^{2}-(1-\alpha_{1})(1-\alpha_{2})}.       \label{eq:eig12}
\end{align}

\begin{rem}\label{rem:a1a2}%
    The property $p\leq q$ was used to arrive at these results. If instead $p > q$,
    we reverse the definitions of $\Proj_\uset$ and $\Proj_\vset$ in Proposition~\ref{prp:projections}.
    Noting that $\spectrum{T}=\spectrum{T^\top}$,
    we get a new block-diagonal matrix $\bar T$ with blocks
    $\bar{T}_1=T_1^\top$, $\bar{T}_3=T_3^\top$ and $\bar{T}_{2} = (1-\alpha_{1})I_{p-q}$.
    Therefore, the matrix can have eigenvalues
    $1-\alpha_1$ or $1-\alpha_2$ depending on the dimensions of $\uset$ and $\vset$.
\end{rem}

If either $p=0$ or $q=0$, then the problem is trivial.
We note that if $p+q\geq n$, we can simply embed the sets in a bigger space. Since $\uset$ and $\vset$ are contained in the original space, the iterates will also stay in this subspace if the initial point is. The algorithm therefore behaves identically and the extra dimensions can be ignored. Although we do not have an explicit expression for the operator $T$ in this case, we can calculate the eigenvalues, as stated in the following theorem.

\begin{thm}\label{thm:gap-all-eigs}
    Let $\uset$ and $\vset$ be subspaces in $\realR^{n}$
    and let $p\ldef \dim(\uset)$, $q\ldef \dim(\vset)$,
    and $s\ldef\dim(\uset\cap\vset)$.
    The eigenvalues of
    $T = \Proj_{\uset}^{\alpha_{2}}\Proj_{\vset}^{\alpha_{1}}$ are
    \begin{align*}
        & \{1\}^s, \{(1-\alpha_1)(1-\alpha_2)\}^{s+n-p-q},\\
        & \{1-\alpha_2\}^{\max(0,q-p)}, \{1-\alpha_1\}^{\max(0,p-q)},\\
        & \{\lambda_{i}^{1,2}\}\, \text{for every $i\in\{s+1,\ldots,\min(p,q)\}$ }
    \end{align*}
    where $\lambda_{i}^{1,2}$ is defined by \eqref{eq:eig12} and
    $\{\lambda\}^i$ denotes (possibly zero) multiplicity $i$ of eigenvalue $\lambda$.
\end{thm}
\begin{pf}
    When either $p=0$ or $q=0$, we get $s=0$ and the result is trivial from the definition of the projections and $T$.
    The case when $p\leq q$ and $p+q<n$ follows directly from Proposition \ref{prp:projections} by observing that $s$ of the eigenvalues in $1$ and $(1-\alpha_1)(1-\alpha_2)$ arise from $\lambda_i^{1,2}$ for $i\in\{1,\dots,s\}$, i.e., when $\theta_i=0$.

    For the case when $q<p$ and $p+q<n$ it follows from Remark \ref{rem:a1a2} that the eigenvalues will be $1-\alpha_1$ instead of $1-\alpha_2$, and that the rest of the eigenvalues are the same.

    For the case when $p+q\geq n$ we provide a proof similar to that in~\cite[p. 54]{Bauschke_lin_rate_Friedrich}.
    We can extend the space $\realR^n$ to $\realR^{n+k}\ldef\realR^n\times\realR^k$ so that $p+q < n +k \rdef \bar n$,
    where we define the scalar product in this new space as $\langle(u_1,u_2),(v_1,v_2)\rangle\ldef\langle u_1,v_1\rangle+\langle u_2,v_2\rangle$ for $u_1,v_1\in\realR^n, u_2,v_2\in\realR^k$.

    Let $\bar{\uset}\ldef\uset\times\{0_k\}$, $\bar{\vset}\ldef\vset\times\{0_k\}$ so that
    \begin{align*}
    \Proj_{\bar\uset} =
        \begin{pmatrix}
    \Proj_{\uset} & 0 \\
    0 & 0_k
    \end{pmatrix},\quad
    \Proj_{\bar\vset} =
        \begin{pmatrix}
    \Proj_{\vset} & 0 \\
    0 & 0_k
    \end{pmatrix}.
    \end{align*}
    It follows that
    \begin{align}\label{eq:textended}
    \bar T \ldef \Proj_{\bar\uset}^{\alpha_2}\Proj_{\bar\vset}^{\alpha_1} = \begin{pmatrix}
    T & 0 \\
    0 & (1-\alpha_1)(1-\alpha_2)I_k
    \end{pmatrix},
    \end{align}
    where $T=\Proj_{\uset}^{\alpha_2}\Proj_{\vset}^{\alpha_1}$. $\bar{T}$ has the same eigenvalues as $T$, as well as $k$ new eigenvalues in $(1-\alpha_1)(1-\alpha_2)$.
    As seen in the definition of $\bar{\uset},\bar{\vset}$ and $\bar{T}$, these \emph{artificial} eigenvalues correspond to directions that are orthogonal to the original space $\realR^n$.
    If we now apply the result for $p+q<\bar{n}$ to $\bar{T}$, and observe that the principal angles are the same for $\bar{\uset},\bar{\vset}$ as for $\uset,\vset$, we see that the eigenvalues are as those stated in the theorem, but with $s+\bar{n}-p-q$ eigenvalues in $(1-\alpha_1)(1-\alpha_2)$. Subtracting the $k$ \emph{artificial} eigenvalues, we conclude that the operator $T$ must have $s+n-p-q$ eigenvalues in $(1-\alpha_1)(1-\alpha_2)$.
\end{pf}
\begin{prp}\label{prp:singularpq}
    Let $\uset$ and $\vset$ be subspaces in $\realR^{n}$ and let $p= \dim(\uset)$, $q= \dim(\vset)$, and $s=\dim(\uset\cap\vset)$. Then the GAP operator $S$ satisfies
    \begin{align*}
        \sigma(S)&=\|S-S^\infty\|\\
        &\leq \max(\|S_1-S_1^\infty\|,|1-\alpha_2(1-\alpha)|,|\alpha+(1-\alpha)(1-\alpha_1)(1-\alpha_2)|, |1-\alpha|),
    \end{align*}
    where $S_1=(1-\alpha)I+\alpha T_1$ with $T_1$ defined in \eqref{eq:T1}.
\end{prp}
\begin{pf}
    If either $p=0$ or $q=0$ we trivially have $S=(1-\alpha)I$ so $\|S-S^\infty\|=|1-\alpha|$ and the result holds.
    If $p<q$ and $p+q<n$, $p,q\geq1$ then it follows directly from Proposition \ref{prp:projections} with $S_i=(1-\alpha)I + \alpha T_i$ that
    \begin{align*}
        \|S-S^\infty\|&=\|D\left((1-\alpha)I +\alpha T\right)D^*-\left(D((1-\alpha)I +\alpha T)D^*\right)^\infty\|\\
        &=\|((1-\alpha)I +\alpha T) - ((1-\alpha)I +\alpha T)^\infty)\|\\
                      &=\left\|\begin{bmatrix}S_1-S_1^\infty&0&0\\0&S_2-S_2^\infty&0\\0&0&S_3-S_3^\infty\end{bmatrix}\right\|\\
        &\leq\max(\|S_1-S_1^\infty\|,|1-\alpha_2(1-\alpha)|,|\alpha+(1-\alpha)(1-\alpha_1)(1-\alpha_2)|)
    \end{align*}
    and the result holds.
    If $p<q$ and $p+q \geq n$ we extend the space as in the proof of Theorem \ref{thm:gap-all-eigs}. Since $\bar{T}$ in \eqref{eq:textended} is a block diagonal matrix containing $T$ we get with $\bar{S}=(1-\alpha)I+\alpha\bar{T}$ that $\|S-S^\infty\|\leq\|\bar{S}-\bar{S}^\infty\|$ and the result follows by applying the case $p+q<n$ to the operator $\bar{S}$.
    For the remaining cases where $p<q$, we note as in Remark \ref{rem:a1a2} that we can study $S^\top=(1-\alpha)I+\alpha\Proj_{\vset}^{\alpha_{1}}\Proj_{\uset}^{\alpha_{2}}$ where the relative dimensions of the subspaces now satisfy the assumptions. Applying the previous results to this case yields $\|S^\top-{S^\top}^\infty\|=\|(S-S^\infty)^\top\|=\|S-S^\infty\|$ and the proof is complete.
\end{pf}
It was shown in \cite{falt-optimal} that the parameters
\begin{align}\label{eq:optpar}
    \alpha=1,\quad\alpha_{1}=\alpha_{2}=\alpha^{*}\ldef \frac{2}{1+\sin{\theta_{F}}},
\end{align}
result in that the subdominant eigenvalues of $S$ have magnitude $\gamma(S)=\gamma^{*}$, where
\begin{align}\label{eq:optrate}
\gamma^* \ldef \alpha^*-1 = \frac{1-\sin(\theta_F)}{1+\sin(\theta_F)}.
\end{align}
When the Friedrichs angle does not exist, i.e., when one subspace is contained in the other,
we define $\alpha^*=1$ and $\gamma^*=0$.
\mnoteh{Here, and in Theorem 2, we need that $\gamma^*$ is defined when there is no $\theta_F$..}
The next two theorems show that this rate is optimal under mild assumptions.
The theorems were published without proofs by the authors in \cite{falt-optimal}.
We restate them with minor modifications and prove them here.

\begin{thm}\label{thm:gapP1-1}\cite[Thm. 1]{falt-optimal}
The GAP operator $S$ in~\eqref{eq:GAP2} for the feasibility problem $(\uset,\vset)$ with linear subspaces $\uset,\vset\subseteq\realR^n$ and
with $\alpha,\alpha_1,\alpha_2$ as defined in~\eqref{eq:optpar}
satisfies $\gamma(S)=\gamma^{*}$, where $\gamma(S)$ and $\gamma^*$ are defined in Definition~\ref{def:subdominant} and \eqref{eq:optrate} respectively.
Moreover, $S$ is linearly convergent with any rate $\mu\in\left(\gamma^{*},1\right)$.
\end{thm}
\begin{pf}
  See Appendix~\ref{app:gapP1-1}.
\end{pf}

\begin{rem}\label{rem:larger-angle}
    Although the rate in Theorem \ref{thm:gapP1-1} is dependent on knowing the true Friedrichs angle $\theta_F$, it is sufficient to have some conservative estimate $\hat{\theta}_F<\theta_F$. As seen in the proof of Theorem \ref{thm:gapP1-1}, choosing the parameters as $\alpha_1=\alpha_2=2/(1+\sin\hat{\theta}_F)$, results in the rate $\gamma=(1-\sin\hat{\theta}_F)/(1+\sin\hat{\theta}_F)$.
\end{rem}

Under the assumption that the relative dimensions of the subspaces are unknown,
it was stated that the rate $\gamma^{*}$ is optimal. We restate it with slight modifications for clarity, and prove it here.

\begin{thm}\label{thm:gapiff}\cite[Thm. 2]{falt-optimal}
Let $(\uset_1,\vset_1)$ and $(\uset_2,\vset_2)$ be two feasibility problems,
where the sets are linear subspaces in $\realR^n$.
Assume that $\dim(\uset_1)<\dim(\vset_1)$, $\dim(\uset_2)>\dim(\vset_2)$ and that $c(\uset_1,\vset_1) = c(\uset_2,\vset_2)$
$=\cos(\theta_F)$,
$\theta_F<\pi/2$.
Let $S_1,S_2$ be the corresponding GAP operators as defined in \eqref{eq:GAP2}, both defined with the same parameters $\alpha_1,\alpha_2,\alpha>0$.
Then, both $S_1$ and $S_2$ are linearly convergent with all rates $\mu\in(\gamma^*,1)$\fnote{, for $\gamma^*$ in \eqref{eq:optrate},}
if and only if
\begin{equation*}
    \alpha=1,\quad\alpha_{1}=\alpha_{2}=\alpha^{*}\ldef \frac{2}{1+\sin{\theta_{F}}}.
\end{equation*}
\end{thm}
\begin{pf}
    See Appendix~\ref{app:gapiff}.
\end{pf}

This theorem shows that there is no choice of parameters that can perform better than that in \eqref{eq:optpar} independently of the dimensions of the subspaces.
Any choice of parameters that performs better than those in \eqref{eq:optpar} for a specific problem, where the dimensions of the subspaces are not the same,
will necessarily perform worse on all problems where the relative dimensions are reversed, if the Friedrichs angle is kept constant.

\begin{rem}
    The are a few cases that are excluded in the theorem that should be explained.
    When $\theta_F=\pi/2$, we have $\gamma^*=0$, which is obviously optimal, however, there are choices of $\alpha,\alpha_1,\alpha_2$ other than \eqref{eq:optpar} that achieve this rate.
    The same is true if the Friedrichs angle is not well defined, i.e., when one set is contained in the other.
    In that case, by defining $\theta_F=\pi/2$,
    we get $\gamma(S)=0$ with the parameters in \eqref{eq:optpar},
    but the solution is not unique.

    As noted in \cite{falt-optimal}, there are specific choices of $(\uset,\vset)$ where it is possible to get $\gamma(S)<\gamma^*$.
    However, if one of the principal angles is large enough, for example $\theta_i=\pi/2$,
    then it is not possible to get a rate better than $\gamma^*$.
    In the cases where $\gamma(S)<\gamma^*$,
    the difference in rate is negligible if $\theta_F$ is small, as long as the parameters are chosen so that the algorithm is convergent for every $(\uset,\vset)$.
    For example, if $\dim\uset\leq\dim\vset$ and \emph{all} principal angles $\theta_i$ are small enough,
    then the parameter choice \emph{GAP$2\alpha$} in \cite{falt-optimal}
    \begin{equation*}
        \alpha = 1, \quad \alpha_1=2, \quad \alpha_2=\frac{2}{1+\sin(2\theta_F)}
    \end{equation*}
    achieves a rate of
    \[
        \frac{\cos(\theta_F)-\sin(\theta_F)}{\cos(\theta_F)+\sin(\theta_F)}=1-2\theta_F+2\theta_F^2-8\theta_F^3/3+O(\theta_F^4)\quad (\text{as } \theta_F\rightarrow 0)
    \]
    compared to
    \begin{equation*}
        \gamma^*=\frac{1-\sin(\theta_F)}{1+\sin(\theta_F)}=1-2\theta_F+2\theta_F^2-5\theta_F^3/3+O(\theta^4)\quad (\text{as } \theta_F\rightarrow 0).
    \end{equation*}
    This should be contrasted to the rates of alternating projections and Douglas--Rachford, which are $1-\theta_F^2+O(\theta_F^4)$ and $1-\theta_F^2/2+O(\theta_F^4)$ as $\theta_F\rightarrow 0$ respectively. So for small angles $\theta_F$, the improvement over AP and DR is significant ($O(\theta_F)$), and the difference to \emph{GAP$2\alpha$} is very small ($O(\theta_F^3)$). As mentioned above, the rate for \emph{GAP$2\alpha$} is only valid under an assumption on the relative dimensions of the manifolds, and that all principal angles are small enough.
\end{rem}

\section{Manifolds}
In this section, we study the local properties of the GAP operator applied to two manifolds $\M$ and $\N$ instead of linear subspaces. These results generalize the results in Section 4 of \cite{lewis-2008}, from alternating projections to the GAP algorithm,
with similar proofs but under the relaxed Assumption \ref{ass:regularity} instead of transversality.

We begin by showing that the GAP operator is locally well defined and well behaved around all points that satisfy the regularity assumptions.

\begin{lem}\label{lem:well-behaved}%
Let $\M$ and $\N$ be manifolds that satisfy Assumption~\ref{ass:regularity} at $\bar{x}\in\M\cap\N$, let $\alpha_1,\alpha_2\in[0,2]$, and let $\alpha\in\reals$.
Then $\Proj_{\M\cap\N}$, $\Proj^{\alpha_2}_{\M}\Proj^{\alpha_1}_{\N}$, and $S=(1-\alpha)I+\alpha\Proj^{\alpha_2}_{\M}\Proj^{\alpha_1}_{N}$ are well defined and of class $\C^{k-1}$ around $\bar{x}$.
\end{lem}
\begin{pf}
    From Assumption \ref{ass:regularity} \ref{ass:reg-smooth} it follows that $\M\cap\N$ is a $\C^{k}$  manifold (with $k\geq2$) so from Lemma \ref{lem:proj-manifold-relaxed} we know that there exists $\delta>0$
so that $\Proj_{\M}$, $\Proj_{\N}$, and $\Proj_{\M\cap \N}$ are well defined and of class $\C^{k-1}$
on $\Ball_{\delta}(\bar{x}).$ Let $x\in \Ball_{\delta/3}(\bar{x})$. Then
\begin{align*}
\left\Vert \bar{x}-\Proj_{\N}^{\alpha_1}(x)\right\Vert &\leq\left\Vert \bar{x}-x\right\Vert +\left\Vert x-\Proj_{\N}^{\alpha_1}(x)\right\Vert = \left\Vert \bar{x}-x\right\Vert +\alpha_1\left\Vert x-\Proj_{\N}(x)\right\Vert \\
& \leq\left\Vert \bar{x}-x\right\Vert +\alpha_1\left\Vert x-\bar{x}\right\Vert\leq 3\left\Vert x-\bar{x}\right\Vert \leq\delta
\end{align*}
so $\Proj_{\N}^{\alpha_1}(x)\in \Ball_{\delta}(\bar{x})$ and therefore 
$\Proj^{\alpha_2}_{\M}\Proj^{\alpha_1}_{N}$ and $S$ are well defined and $\C^{k-1}$ on $\Ball_{\delta/3}(\bar{x})$.
\end{pf}

To simplify notation, we denote the GAP operator applied to the tangent spaces $\T_{\M}(\bar{x})$ and $T_{\N}(\bar{x})$ by
\begin{align}\label{eq:gap-defin-tangent}
S_{\T(\bar{x})}\ldef (1-\alpha)I+\alpha\Proj_{\T_{\M}(\bar{x})}^{\alpha_{2}}\Proj_{\T_{\N}(\bar{x})}^{\alpha_{1}}.
\end{align}
We next show that the local behavior of $S$ around a point $\bar{x}\in\M\cap\N$ can be described by $S_{\T(\bar{x})}$.
\begin{lem}\label{lem:gap-tangent}%
    Let $\M$ and $\N$ be manifolds that satisfy Assumption \ref{ass:regularity} at $\bar{x}\in\M\cap\N$. Then the Jacobian at $\bar{x}\in\reals^n$ of the GAP operator $S$ in~\eqref{eq:GAP2} is given by
    \[
      \jacobian{S}{\bar{x}}=S_{\T(\bar{x})},
    \]
    with $S_{\T(\bar{x})}$ defined in \eqref{eq:gap-defin-tangent}.
\end{lem}
\begin{pf}
    Using the chain rule, $\bar{x}\in\M\cap\N$, and Lemma \ref{lem:proj-manifold-relaxed} we conclude
\begin{align*}
\jacobian{\Proj^{\alpha_2}_\M\Proj^{\alpha_1}_\N}{\bar{x}}&=\jacobian{\Proj^{\alpha_2}_\M}{\Proj^{\alpha_1}_\N(\bar{x})}\jacobian{\Proj^{\alpha_1}_\N}{\bar{x}}
=\jacobian{\Proj^{\alpha_2}_\M}{\bar{x}}\jacobian{\Proj^{\alpha_1}_\N}{\bar{x}}=\Proj_{\T_{\M}(\bar{x})}^{\alpha_{2}}\Proj_{\T_{\N}(\bar{x})}^{\alpha_{1}}.
\end{align*}
Moreover
\begin{align*}
    \jacobian{S}{\bar{x}}&=\jacobian{(1-\alpha)I}{\bar{x}}+\alpha\jacobian{\Proj^{\alpha_2}_\M\Proj^{\alpha_1}_\N}{\bar{x}}=(1-\alpha)I+\alpha\Proj_{\T_{\M}(\bar{x})}^{\alpha_{2}}\Proj_{\T_{\N}(\bar{x})}^{\alpha_{1}}=S_{\T(\bar{x})}
\end{align*}
by definition of $S_{\T(\bar{x})}$ in \eqref{eq:gap-defin-tangent}.
\end{pf}

%
\begin{prp}\label{prp:intersection-limit}
    Let $\M$ and $\N$ be manifolds that satisfy Assumption \ref{ass:regularity} at $\bar{x}\in\M\cap\N$ and let the parameters of the GAP operator $S$ in \eqref{eq:GAP2} satisfy Assumption \ref{ass:alpha} case \ref{ass:alpha-1} or \ref{ass:alpha-2}. Then
    \begin{align}
      \T_{\M(\bar{x})\cap \N(\bar{x})} = \T_{\M(\bar{x})}\cap T_{\N(\bar{x})} = \fix S_{\T(\bar{x})}
      \label{eq:tangent-space_fixed-points}
    \end{align}
    and
    \begin{align}
      \Proj_{\fix S_{\T(\bar{x})}} = S_{\T(\bar{x})}^\infty.
      \label{eq:projection_limit-point}
    \end{align}
\end{prp}
\begin{pf}
    The first equality follows from Assumption \ref{ass:regularity}. From Lemma \ref{lem:fixed-points}, under Assumption \ref{ass:alpha} cases \ref{ass:alpha-1} and \ref{ass:alpha-2}, we know that $\fix S_{\T(\bar{x})}=\T_{\M(\bar{x})}\cap T_{\N(\bar{x})}$ and \eqref{eq:tangent-space_fixed-points} is proven. From Theorem~\ref{thm:gapP1-1}, we know that $S_{\T(\bar{x})}$ is convergent (to $S_{\T(\bar{x})}^{\infty}$)  and non-expansiveness of $S_{\T(\bar{x})}$ and \cite[Corollary 2.7]{Bauschke_opt_rate_matr} imply that
    $\Proj_{\text{Fix}S_\T(\bar{x})}=S_{\T(\bar{x})}^{\infty}$ and \eqref{eq:projection_limit-point} is proven.
\end{pf}
We next show that the convergence rate of $S^k(x)$ to the intersection tends to the rate $\gamma(S_{\T(\bar{x})})$ as the initial point gets closer to the intersection and the number of iterations $k$ increases.

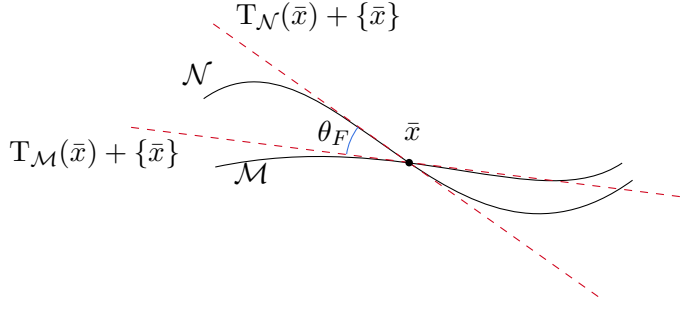
\begin{figure}
  \begin{center}
  \begin{tikzpicture}[x=0.75pt,y=0.75pt,yscale=-1,xscale=1]

  \draw    (198.67,68.33) .. controls (261.67,23.33) and (322.67,175.33) .. (412.67,109.33) ;
  \draw    (407.47,100.7) .. controls (367.47,127.7) and (298.47,82.7) .. (204.47,102.7) ;
  \draw [color={rgb, 255:red, 208; green, 2; blue, 27 }  ,draw opacity=1 ] [dash pattern={on 3pt off 3pt}]  (203.47,30.7) -- (397.62,169.4) ;
  \draw [color={rgb, 255:red, 208; green, 2; blue, 27 }  ,draw opacity=1 ] [dash pattern={on 3pt off 3pt}]  (162.46,82.74) -- (438.47,117.86) ;
  \draw  [color={rgb, 255:red, 54; green, 114; blue, 226 }  ,draw opacity=1 ] (270,96.23) .. controls (270.55,91.77) and (272.61,86.73) .. (275.78,82.53) ;

  \draw (273.00,85.54) node [anchor=east] [inner sep=0.75pt]    {$\theta_{F}$};
  \draw (212.42,99.58) node [anchor=north west][inner sep=0.75pt]    {$\M$};
  \draw (297.42,79.58) node [anchor=north west][inner sep=0.75pt]    {$\bar{x}$};
  \draw (187.42,49.58) node [anchor=north west][inner sep=0.75pt]    {$\N$};
  \draw (213.42,18.58) node [anchor=north west][inner sep=0.75pt]    {$\T_\N(\bar{x})+\{\bar{x}\}$};
  \draw (100.42,87.58) node [anchor=north west][inner sep=0.75pt]    {$\T_\M(\bar{x})+\{\bar{x}\}$};

  \node[fill=black,circle,outer sep=0, inner sep=1] at (301.19, 100.5) {};
  \end{tikzpicture}
  \caption{Illustration of manifolds $\M$ and $\N$ in $\reals^2$ and the approximation by tangent planes at a point $\bar{x}\in\M\cap\N$.}
  \end{center}
\end{figure}

\begin{thm}\label{thm:limitrate1}
Let $\M$ and $\N$ be manifolds that satisfy Assumption \ref{ass:regularity} at $\bar{x}\in\M\cap\N$ and let the parameters of the GAP operator $S$ in \eqref{eq:GAP2} satisfy Assumption \ref{ass:alpha} case \ref{ass:alpha-1} or \ref{ass:alpha-2}.
Then
\begin{enumerate}[label=\arabic*.,ref=\arabic*]
    \item\label{thm:limitrate1:case1} for all
      $c>\left\Vert S_{\T(\bar{x})}-\Proj_{\T_{\M}(\bar{x})\cap T_{\N}(\bar{x})}\right\Vert$ with $S_{\T(\bar{x})}$ is defined in \eqref{eq:gap-defin-tangent},
        there exists some $\eta>0$ so that for all $x\in\Ball_\eta(\bar{x})$
        \begin{equation}\label{eq:contraction-rate}
        \left\Vert S(x)-\Proj_{\M\cap \N}(x)\right\Vert \leq c\left\Vert x-\Proj_{\M\cap \N}(x)\right\Vert.
        \end{equation}
    \item\label{thm:limitrate1:case2} for all
        $\mu_{\bar{x}}\in(\gamma(S_{\T(\bar{x})}),1)$ there exists $N\in\naturalN$, such that for any $k\geq N$
        \begin{equation}
        \limsup_{x\rightarrow\bar{x},x\not\in\M\cap\N}\frac{\left\Vert S^{k}(x)-\Proj_{\M\cap \N}(x)\right\Vert }{\left\Vert x-\Proj_{\M\cap \N}(x)\right\Vert }\leq\mu_{\bar{x}}^{k}.\label{eq:limk}
        \end{equation}

\end{enumerate}
\end{thm}
\begin{pf}
  Let $x_{r}\not\in\M\cap \N$ and denote $\bar{x}_r=\Proj_{\M\cap \N}(x_r)$.
Since $\bar{x}_{r}\in\M\cap\N$ we trivially have $S\bar{x}_r=\bar{x}_r$.

Let $x_r$ be in the region around $\bar{x}$ for which $S$ and $\Proj_{\M\cap \N}$ are well-defined and $\mathcal{C}^1$ according to Lemma~\ref{lem:well-behaved}.
By \cite[Eq (3.8.1), Thm 3.8.1]{cartan71}, a $\mathcal{C}^1$ function $f:\realR^n\rightarrow\realR^n$ at a point $a\in\realR^n$ can be approximated as
\begin{align*}
f(x)-f(y) = \jacobian{f}{a}(x-y) + \|x-y\|\psi(x,y), \text{ where } \lim_{x,y\rightarrow a}\psi(x,y)=0,
\end{align*}
at $x,y\in\realR^n$.
Using this, with $f(x)=S(x)-\Proj_{\M\cap \N}(x)$, at $x=x_r, y=\bar{x}_r, a=\bar{x}$ we get
\begin{align}\label{eq:expansion-S}
S(x_{r})-\Proj_{\M\cap \N}(x_{r})&=(\jacobian{S}{\bar{x}} - \jacobian{\Proj_{\M\cap \N}}{\bar{x}})(x_{r}-\bar{x}_r)+ \|x_r-\bar{x}_r\|\psi(x_r,\bar{x}_r),\\ &\text{ where } \lim_{x_r,\bar{x}_r\rightarrow \bar{x}}\psi(x_r,\bar{x}_r)=0.\nonumber
\end{align}
We can replace the Jacobians by noting that Lemma \ref{lem:gap-tangent}, Lemma \ref{lem:proj-manifold}, and Assumption \ref{ass:regularity} \ref{ass:reg-intersection} at $\bar{x}\in\reals^n$ imply
\begin{align*}
    \jacobian{S}{\bar{x}} - \jacobian{\Proj_{\M\cap \N}}{\bar{x}} = S_{\T(\bar{x})} - \Proj_{\T_{\M}(\bar{x})\cap T_{\N}(\bar{x})}.
\end{align*}
Using this equality in \eqref{eq:expansion-S}, taking the norm of both sides, applying the triangle inequality and Cauchy--Schwarz, and dividing by $\|x_r-\bar{x}_r\|$ result in
\begin{equation}
\frac{\left\Vert S(x_{r})-\bar{x}_{r}\right\Vert }{\left\Vert x_{r}-\bar{x}_{r}\right\Vert }\leq\left\Vert S_{\T(\bar{x})}-\Proj_{\T_{\M}(\bar{x})\cap T_{\N}(\bar{x})}\right\Vert +\Vert\psi(x_r,\bar{x}_r)\Vert, \text{ if }x_r\neq\bar{x}_r.
\end{equation}
Continuity of $\Proj_{\M\cap \N}$ around $\bar{x}$ means that $\psi(x_r,\bar{x}_r)=\psi(x_r,\Proj_{\M\cap \N}(x_r))\rightarrow0$ as $x_r\rightarrow \bar{x}$,
so for any $c>\left\Vert S_{\T(\bar{x})}-\Proj_{\T_{\M}(\bar{x})\cap T_{\N}(\bar{x})}\right\Vert$,
there exists some $\eta>0$ so that
\begin{equation}
\forall x_r\in\Ball_\eta(\bar{x}):\quad
\left\Vert S(x_{r})-\bar{x}_{r}\right\Vert \leq c\left\Vert x_{r}-\bar{x}_{r}\right\Vert.
\end{equation}
This proves part \ref{thm:limitrate1:case1} of the theorem.

In the same way for $S^k$, since $S(\bar{x})=S_{\T(\bar{x})}(\bar{x})=\bar{x}$, using the chain rule, we get
\[
\jacobian{S^{k}}{\bar{x}}=\left(\jacobian{S}{\bar{x}}\right)^{k}=S_{\T(\bar{x})}^{k},
\]
so in the same way we conclude
\begin{equation}
\frac{\left\Vert S^k(x_{r})-\bar{x}_{r}\right\Vert }{\left\Vert x_{r}-\bar{x}_{r}\right\Vert }\leq\left\Vert S^k_{\T(\bar{x})}-\Proj_{\T_{\M}(\bar{x})\cap T_{\N}(\bar{x})}\right\Vert +\psi(x_r,\bar{x}_r), \text{ if }x_r\neq\bar{x}_r
\end{equation}
From Proposition \ref{prp:intersection-limit} we have that
$\Proj_{\T_{\M}(\bar{x})\cap T_{\N}(\bar{x})}=S_{\T(\bar{x})}^{\infty}$ and thus
\begin{align*}
\frac{\left\Vert S^{k}(x_{r})-\bar{x}_{r}\right\Vert }{\left\Vert x_{r}-\bar{x}_{r}\right\Vert } \leq \left\Vert S_{\T(\bar{x})}^{k}-S^\infty_{\T(\bar{x})}\right\Vert  +\psi(x_r,\bar{x}_r), \text{ if }x_r\neq\bar{x}_r.
\end{align*}
Continuity of $\Proj_{\M \cap \N}$ around $\bar{x}=\Proj_{\M \cap \N}(\bar{x})$, with $\bar{x}_r=\Proj_{\M \cap \N}(x_r)$, implies
\begin{align*}
    \limsup_{x_r\rightarrow\bar{x},x_r\not\in\M\cap\N} \frac{\left\Vert S^{k}(x_{r})-\bar{x}_{r}\right\Vert }{\left\Vert x_{r}-\bar{x}_{r}\right\Vert } & \leq
    \left\Vert S_{\T(\bar{x})}^{k}-S^\infty_{\T(\bar{x})} \right\Vert.
\end{align*}
Using the results in \cite{falt-optimal} with Definitions \ref{def:linear}, \ref{def:limit}, \ref{def:subdominant}, and Facts \ref{fct:limitexists}, \ref{fact:Convergent-at-rate}
imply that for any $\mu_{\bar{x}}$ with $\gamma(S_{\T(\bar{x})})<\mu_{\bar{x}}$
there exists $N\in\naturalN$ so that for all $k\geq N$
\[
\left\Vert S_{\T(\bar{x})}^{k}-S_{\T(\bar{x})}^{\infty}\right\Vert \leq\mu_{\bar{x}}^{k}.
\]
We conclude that for any $\mu_{\bar{x}}\in(\gamma(S_{\T(\bar{x})}),1)$, there exists $N$ such that for all $k\geq N$
\begin{equation}
\limsup_{x\rightarrow\bar{x},x\not\in \M\cap \N}\frac{\left\Vert S^{k}(x)-\Proj_{\M\cap\N}(x)\right\Vert }{\left\Vert x-\Proj_{\M\cap\N}(x)\right\Vert }\leq\mu_{\bar{x}}^{k},
\end{equation}
which proofs part \ref{thm:limitrate1:case2} of the theorem.
%
\end{pf}

It remains to show that the sequence of iterates actually converges. To do this, we first show that $\|S_{\T(\bar{x})}-\Proj_{\T_\M(\bar{x})\cap \T_\N(\bar{x})}\|<1$.
\begin{lem}\label{lem:contraction}%
Let $\alpha,\alpha_1,\alpha_2$ satisfy Assumption \ref{ass:alpha} case \ref{ass:alpha-1} or \ref{ass:alpha-2}, and let $\M$ and $\N$ be manifolds that satisfy Assumption \ref{ass:regularity} at $\bar{x}\in\M\cap\N$.
Then
\begin{equation}
    \rate{S_{\T(\bar{x})}}\ldef \|S_{\T(\bar{x})}-\Proj_{\T_\M(\bar{x})\cap \T_\N(\bar{x})}\|<1
\end{equation}
where $S_{\T(\bar{x})}=\alpha \Proj^{\alpha_2}_{\T_\M(\bar{x})}\Proj^{\alpha_1}_{\T_\N(\bar{x})} + (1-\alpha)I$
\end{lem}
\begin{pf}
    First note that $\Proj_{\T_\M(\bar{x})\cap \T_\N(\bar{x})}=\Proj_{\text{Fix}S_{\T(\bar{x})}}=S_{\T(\bar{x})}^{\infty}$ by Proposition \ref{prp:intersection-limit}.
    Proposition \ref{prp:singularpq} therefore gives that
    \begin{align*}
        \|S_{\T(x)}-S_{\T(x)}^\infty\|\leq \max(&\|S_1-S_1^\infty\|,|1-\alpha_2(1-\alpha)|,|\alpha+(1-\alpha)(1-\alpha_1)(1-\alpha_2)|, |1-\alpha|),
    \end{align*}
    where $S_1$ is a block diagonal matrix with blocks $S_{1_i}=(1-\alpha)I+\alpha T_{1_i}$ and $T_{1_i}$ are defined in \eqref{eq:T1matrix}. 
    Under Assumption \ref{ass:alpha} case \ref{ass:alpha-1} or \ref{ass:alpha-2}
    we have $|1-\alpha_2(1-\alpha)|<1$, $|\alpha+(1-\alpha)(1-\alpha_1)(1-\alpha_2)|<1$ and $|1-\alpha|< 1$. It remains to show that $\|S_1-S_1^\infty\|=\max_i{\|S_{1_i}-S_{1_i}^\infty\|}<1$.
    We now consider each block $S_{1_i}$ corresponding to each of the principal angles $\theta_i$.
    Each block with $\theta_i=0$ becomes
    \begin{align*}
        S_{1_i}&=\alpha T_{1_i}+(1-\alpha)I =
        \begin{pmatrix}
        1 & 0 \\
        0 & \alpha(1-\alpha_1)(1-\alpha_2)+(1-\alpha)
        \end{pmatrix}\\
        S_{1_i}^\infty &= \begin{pmatrix}
       1 & 0 \\
       0 & 0
       \end{pmatrix},
   \end{align*}
   so the corresponding singular values are $0$ and $|\alpha(1-\alpha_1)(1-\alpha_2)+(1-\alpha)|<1$.
   The remaining cases are $\theta_i\in(0,\pi/2]$ for which $(S_{1_i})^\infty=\Proj_{\fix S_{1_i}} = 0$.
   To study the largest singular value $\|S_{1_i}-S_{1_i}^\infty\|=\|S_{1_i}\|=\|\alpha T_{1_i} + (1-\alpha)I\|$ \fnote{we first note that $S_{1_i}$ is non-expansive} so
   $\|S_{1_i}\|\leq 1$, hence we only need to show that $\|S_{1_i}\|\neq 1$.
   From the triangle inequality we get $\|\alpha T_{1_i} + (1-\alpha)I\|\leq \alpha\| T_{1_i}\|+(1-\alpha)\leq 1$,
   with equality only if $\|T_{1_i}\|=1$. To this end, we consider
   $\|T_{1_i}\|^2=\rho(T_{1_i}T_{1_i}^\top)$ 
   and study the eigenvalues of $T_{1_i}T_{1_i}^\top$.
    Non-expansiveness again implies that $\|T_{1_i}\|\leq 1$.
    We now aim to show that these blocks have singular values smaller than $1$ when $\theta_i\in(0,\pi/2]$.
    After simplification using the identity $s_i^2+c_i^2=1$, we get
    \begin{align*}
        T_{1_i}T_{1_i}^\top & =
        \begin{pmatrix}
            1-2\alpha_1s_i^2+\alpha_1^2s_i^2 & (2-\alpha_1) \alpha_1 (1-\alpha_2) c_i s_i \\
             (2-\alpha_1) \alpha_1 (1-\alpha_2) c_i s_i & (1-\alpha_2)^2 (1-2\alpha_1c_i^2+\alpha_1^2c_i^2)
        \end{pmatrix} =:
        \begin{pmatrix}
            a & b \\
            c & d
        \end{pmatrix}.
    \end{align*}
    For any of these eigenvalues to be $1$, it must hold that
    \begin{align*}
        \det\begin{pmatrix}
            a-1 & b \\
            c & d-1
        \end{pmatrix}=0,
    \end{align*}
    i.e.,
    \begin{align}\label{eq:eigval1}
        0 &= 1 -a -d +ad -bc.
    \end{align}
    Simplifying the expressions yields the following identities
    \begin{align*}
        1-a-d &= \alpha_1s_i^2(2-\alpha_1)-(1-\alpha_2)^2(1-2\alpha_1c_i^2+\alpha_1^2c_i^2)\\
        ad &= (1-\alpha_2)^2(\alpha_1^2c_i^2s_i^2(4-4\alpha_1+\alpha_1^2)+(1-\alpha_1)^2)\\
        bc &= (1-\alpha_2)^2\alpha_1^2c_i^2s_i^2(4-4\alpha_1+\alpha_1^2)\\
        ad-bc &= (1-\alpha_1)^2(1-\alpha_2)^2
    \end{align*}
    and thus
    \begin{align*}
        1-a-d+ad-bc &= \alpha_1s_i^2(2-\alpha_1)-(1-\alpha_2)^2(1-2\alpha_1c_i^2+\alpha_1^2c_i^2)\\
        &\quad + (1-\alpha_1)^2(1-\alpha_2)^2\\
        &= s_i^2\alpha_1(2-\alpha_1)-(1-\alpha_2)^2(2\alpha_1(1-c_i^2)+\alpha_1^2(c_i^2-1))\\
        &= s^2\alpha_1(2-\alpha_1)-(1-\alpha_2)^2\alpha_1s_i^2(2-\alpha_1)\\
        &= s_i^2\alpha_1\alpha_2(2-\alpha_1)(2-\alpha_2).
    \end{align*}
    From \eqref{eq:eigval1} we conclude that for the largest eigenvalue to be $1$, it must hold that
    \begin{align*}
        0 &= \sin(\theta_i)^2\alpha_1\alpha_2(2-\alpha_1)(2-\alpha_2).
    \end{align*}
    Within the ranges $\alpha_1,\alpha_2\in(0,2)$ and $\theta_i\in(0,\pi/2]$ we have \begin{align*}
    \sin(\theta_i)^2\alpha_1\alpha_2(2-\alpha_1)(2-\alpha_2)>0,
    \end{align*}
    which leads to 
    $\rho(T_{1_i}T_{1_i}^\top)=\|T_{1_i}\|^2<1$,
    and thus $\|S_{1_i}\|<1$. This completes the proof for case \ref{ass:alpha-1} from Assumption \ref{ass:alpha}.

    Now consider the parts of case \ref{ass:alpha-2} from Assumption \ref{ass:alpha} that are not covered by case \ref{ass:alpha-1}, i.e., $\alpha\in(0,1)$ and $\alpha_1=2$ or $\alpha_2=2$ (but not both) implying that $\|T_{1_i}\|=1$.
    Suppose that $\|S_{1_i}\|=1$.
    From compactness of the unit sphere in $\realR^n$ and continuity of the norm we get from the definition of the operator norm that there exists a $\|v\|=1$ such that
    $\|S_{1_i}v\|=1$. But then $1=\|S_{1_i}v\|^2=\|\alpha T_{1_i}v + (1-\alpha)v\|^2$.
    However, on the boundaries $\alpha=0$ or $\alpha=1$ we get $\|S_{1_i}v\|=1$.
    Since the squared norm is strongly convex we have for any $\alpha\in(0,1)$ where $T_{1_i}v\neq v$ the contradiction $\|\alpha T_{1_i}v + (1-\alpha)v\|^2<1$.
    This leaves the case where $T_{1_i}v = v$, which means that $v$ is a fixed point of $T$, but the only fixed point is $v=0$.
    Thus, there is no $\|v\|=1$ such that $\|S_{1_i}v\|=1$ and therefore $\|S_{1_i}\|<1$. This concludes the proof.
\end{pf}

We are now ready to show that the algorithm will locally converge to some point in the intersection with the contraction factor in Lemma \ref{lem:contraction}.
The proof is similar to that in \cite{lewis-2008}, where the authors show the result for the special case of alternating projections.

\begin{thm}\label{thm:manifold-some-linear-convergence}
    Let $\M$ and $\N$ be manifolds that satisfy Assumption \ref{ass:regularity} at $\bar{x}\in\M\cap\N$, and let the parameters of the GAP operator $S$ in \eqref{eq:GAP2} satisfy Assumption \ref{ass:alpha} case \ref{ass:alpha-1} or \ref{ass:alpha-2}. If the initial
    point $x_{0}\in\reals^n$ is close enough to $\bar{x}$ then the GAP method in Definition~\ref{def:GAP}
    is well defined. Moreover, the sequence $(x_{k})_{k\in\naturalN}$ converges to
    some point $x^{*}\in \M\cap \N$,
    and for every $\mu_{\bar{x}}\in(\rate{S_{\T(\bar{x})}},1)$, there exists a $\beta>0$ such that
    \begin{equation}\label{eq:thm-rate}
    \|x_{k}-x^{*}\|\leq\beta \mu_{\bar{x}}^{k}.
    \end{equation}
\end{thm}
\begin{pf}
    By Lemma \ref{lem:contraction} we have $\rate{S_{\T(\bar{x})}}=\|S_{\T(\bar{x})}-\Proj_{\T_\M(\bar{x})\cap \T_\N(\bar{x})}\|<1$.
    Let $c\in(0,1)$ be such that
    $\|S_{\T(\bar{x})}-\Proj_{\T_\M(\bar{x})\cap \T_\N(\bar{x})}\|<c<1$
    and choose $\eta$ such that $Sx$ and $\Proj_{\M\cap \N}(x)$ are well defined by Lemma~\ref{lem:well-behaved} for $x\in B_{\eta}(\bar{x})$ and so that the conditions of Theorem~\ref{thm:limitrate1} are satisfied. Theorem~\ref{thm:limitrate1} case \ref{thm:limitrate1:case1} then gives
    \begin{equation}\label{eq:ball-contraction}
        \forall x\in B_{\eta}(\bar{x}),\quad\|Sx-\Proj_{\M\cap \N}(x)\|\leq c\|x-\Proj_{\M\cap \N}(x)\|.
    \end{equation}

    Let the initial point $x_{0}\in \Ball_{\delta}(\bar{x})$ where $\delta\ldef\eta/(2\sum_{k=0}^{\infty}c^{k})=\eta(1-c)/2<\eta$ and
    define $\bar{x}_{k}:=\Proj_{\M\cap \N}(x_{k})$.
    By the choice of $\eta$, if $x_{k}\in \Ball_{\eta}(\bar{x})$ then $\bar{x}_{k}$ and $x_{k+1}$ are well defined.
    We now show the following results by induction:
    \begin{align*}
    \|x_{k}-\bar{x}\| & \leq 2\delta\sum_{i=0}^{k}c^{i}\tag{H0}\\
    \|x_{k}-\bar{x}_{k}\| & \leq\delta c^{k}\tag{H1}\\
    \|\bar{x}_{k}-\bar{x}_{k-1}\| & \leq2\delta c^{k}\tag{H2}\\
    \|\bar{x}_{k}-\bar{x}\| & \leq2\delta\sum_{i=0}^{k}c^{i}\tag{H3}
    \end{align*}
    where we note that $2\delta\sum_{i=0}^{k}c^{i}\leq\frac{2\delta}{1-c}=\eta$.

    Case $k=0$: Let $\bar{x}_{-1}\ldef\bar{x}_{0}$. We have trivially
    \begin{align*}
    \|x_{0}-\bar{x}\| & \leq\delta\leq2\delta\tag{H$0^0$}\\
    \|x_{0}-\bar{x}_{0}\| & \leq\|x_{0}-\bar{x}\|\leq\delta\tag{H$1^0$}\\
    \|\bar{x}_{0}-\bar{x}_{-1}\| & =0\leq2\delta\tag{H$2^0$}\\
    \|\bar{x}_{0}-\bar{x}\| & \leq2\delta.\tag{H$3^0$}
    \end{align*}
    Now assume that (H0)-(H3) hold up to some $k$. Then by the triangle inequality, \eqref{eq:ball-contraction},
  (H1), and (H3) we get
    \begin{align*}
    \|x_{k+1}-\bar{x}\|&\leq \|x_{k+1}-\bar{x}_{k}\| + \|\bar{x}_{k}-\bar{x}\|\\
    & \leq
    c\|x_{k}-\bar{x}_{k}\|+\|\bar{x}_{k}-\bar{x}\|
    \leq \delta c^{k+1}+ 2\delta \sum_{i=0}^{k}c^i\leq 2\delta \sum_{i=0}^{k+1}c^i.\tag{H$0^+$}
    \end{align*}
    By the definition of the projection,
    \eqref{eq:ball-contraction}, and (H$1)$ we get
    \begin{align*}
    \|x_{k+1}-\bar{x}_{k+1}\|\leq\|x_{k+1}-\bar{x}_{k}\|\leq c\|x_{k}-\bar{x}_{k}\|\leq\delta c^{k+1}.\tag{H$1^+$}
    \end{align*}
    Again, by the triangle inequality, the definition of projection and (H$1^+$)
    \begin{align*}
    \|\bar{x}_{k+1}-\bar{x}_{k}\|\leq\|\bar{x}_{k+1}-x_{k+1}\|+\|x_{k+1}-\bar{x}_{k}\|\leq2\|x_{k+1}-\bar{x}_{k}\|\leq2\delta c^{k+1}\tag{H$2^+$}
    \end{align*}
    and by (H$2^+)$ and (H$3$):
    \begin{align*}
    \|\bar{x}_{k+1}-\bar{x}\|\leq\|\bar{x}_{k+1}-\bar{x}_{k}\|+\|\bar{x}_{k}-\bar{x}\|\leq2\delta c^{k+1}+2\delta\sum_{i=0}^{k}c^{i}=2\delta\sum_{i=0}^{k+1}c^{i}.\tag{H$3^+$}
    \end{align*}
    By induction we have now shown that (H0)--(H3) must hold for all $k\geq0.$

    We now show that $\left(\bar{x}_{k}\right)_{k\in\naturalN}$ is Cauchy. By the triangle
    inequality, \eqref{eq:ball-contraction}, and (H1):
    \begin{align*}
    \|\bar{x}_{k+1}-\bar{x}_{k}\| & \leq\|\bar{x}_{k+1}-x_{k+1}\|+\|x_{k+1}-\bar{x}_{k}\|\\
     & \leq\|\bar{x}_{k+1}-x_{k+1}\|+c\|x_{k}-\bar{x}_{k}\|\leq\delta c^{k+1}+\delta c^{k+1}\leq2\delta c^{k+1}.
    \end{align*}
    Thus for any $p,k\in\naturalN$ with $p>k$
    \[
    \|\bar{x}_{p}-\bar{x}_{k}\|\leq\sum_{i=k}^{p-1}\|\bar{x}_{i+1}-\bar{x}_{i}\|\leq2\delta\sum_{i=k}^{p-1}c^{i+1}\leq2\delta c^{k+1}\sum_{i=0}^{\infty}c^{i}=\frac{2\delta}{1-c}c^{k+1},
    \]
    so the sequence is Cauchy. Therefore $x^{*}=\lim_{p\rightarrow\infty}\bar{x}_{p}\in\M\cap\N$
    exists and
    \[
    \|x^{*}-\bar{x}_{k}\|\leq\frac{2\delta}{1-c}c^{k+1}.
    \]
    Lastly, by the triangle inequality and (H1)
    \[
    \|x_{k}-x^{*}\|\leq\|x_{k}-\bar{x}_{k}\|+\|\bar{x}_{k}-x^{*}\|\leq\delta c^{k}+\frac{2\delta}{1-c}c^{k+1}=\delta\frac{1+c}{1-c}c^{k},
    \]
    hence \eqref{eq:thm-rate} holds with $\beta=\delta\frac{1+c}{1-c}$ and $\mu_{\bar{x}}=c$.
\end{pf}

Theorem~\ref{thm:manifold-some-linear-convergence} implies that the sequence generated by the generalized alternating projection algorithm converges to a point in the intersection when initiated close enough to a point that satisfies Assumption~\ref{ass:regularity}. However, as is the case for the method of alternating projections, the rate predicted by $\sigma(S_{\T(x^{*})})$ is very conservative. We now show that the iterates converge to the intersection with the faster rate $\gamma(S_{\T(x^{*})})$ from Definition \ref{def:subdominant}. The theorem and proof are similar to that in \cite[Rem. 4]{lewis-2008}, where the authors show it for alternating projections.

\begin{thm}\label{thm:manifold-linear-rate}
    Let $\M$ and $\N$ be manifolds that satisfy Assumption \ref{ass:regularity} at $\bar{x}\in\M\cap\N$, let the initial
    point $x_{0}\in\reals^n$ be close enough to $\bar{x}$, and let the parameters of the GAP operator $S$ in \eqref{eq:GAP2} satisfy Assumption \ref{ass:alpha} case \ref{ass:alpha-1} or \ref{ass:alpha-2}.
    Further assume that $\M$ and $\N$ satisfy Assumption \ref{ass:regularity} at the limit point $x^*$ of the sequence $(x_k)_{k\in\naturalN}$ generated by the GAP method in Definition~\ref{def:GAP}.
    Then the convergence is R-linear to $\M\cap\N$ with any rate $\mu_{x^{*}}\in(\gamma(S_{\T(x^{*})}),1)$.
    That is, for any $\mu_{x^{*}}\in(\gamma(S_{\T(x^{*})}),1)$, there exists $N\in\naturalN$ such that
    \begin{equation}\label{eq:limmanifold}
    d_{\M\cap \N}(x_{k})\leq\mu_{x^{*}}^{k},\quad \forall k > N.
    \end{equation}
\end{thm}
\begin{pf}
    We note that Theorem \ref{thm:manifold-some-linear-convergence} establishes the existence of a limit point $x^*$.
    Take any $\mu_{x^{*}}\in(\gamma(S_{\T(x^{*})}),1)$ and let $\bar{\mu}_{x^{*}} = (\mu_{x^{*}}+\gamma(S_{\T(x^{*})}))/2$.
    Theorem \ref{thm:manifold-some-linear-convergence} implies that eventually $x_{r}\in B_{\eta}(x^{*})$,
    and thus by Theorem~\ref{thm:limitrate1} case \ref{thm:limitrate1:case2}, with $\bar{\mu}_{x^{*}}\in(\gamma(S_{\T(x^{*})}),1)$, there exists
    $N\in\naturalN$ so that $\forall t>N$,
    \begin{align*}
        d_{\M\cap \N}(x_{t+n}) & =\left\Vert S^{t}x_{n}-\Proj_{\M\cap \N}(x_{n})\right\Vert  <\bar{\mu}_{x^{*}}^{t}\left\Vert x_{n}-\Proj_{\M\cap \N}(x_{n})\right\Vert =\bar{\mu}_{x^{*}}^{t}d_{\M\cap \N}(x_{n}),
    \end{align*}
    as long as $x_n\not\in\M\cap\N$.
    By induction this leads to
    \begin{equation}
    d_{\M\cap \N}(x_{kt+n})<\bar{\mu}_{x^{*}}^{kt}d_{\M\cap \N}(x_{n}),\quad\forall k=1,2,3,\dots.\label{eq:distancemanifoldk}
    \end{equation}
    Now fix $t>N$ and assume that \eqref{eq:limmanifold} does not hold,
    then there exists an infinite sequence $r_{1}<r_{2}<\cdots$, all
    satisfying
    \begin{equation}\label{eq:counter-ass}
    d_{\M\cap \N}(x_{r_{j}})>\mu_{x^{*}}^{r_{j}}.
    \end{equation}
    We now show that this is impossible and that the theorem therefore must hold.
    By Lemma \ref{lem:sub-sequence} (see Appendix~\ref{app:lemma}), we can select a sub-sequence $\left(r_{k_{j}}\right)_{j\in\naturalN}$ of $\left(r_{j}\right)_{j\in\naturalN}$
    where we can write $r_{k_{j}}=a+b_{j}t$ for some $a\in\naturalN$ and an increasing
    sequence of integers $\left(b_{j}\right)_{j\in\naturalN}$,
    i.e., we have a new sub-sub-sequence where all iterates are a multiplicity of $t$ iterations apart.
    Thus, picking any
    $b$ so that $a+bt>N$, we have with $r_{k_{j}}=a+b_{j}t=a+bt+(b_{j}-b)t$
    from \eqref{eq:distancemanifoldk} that
    \[
    d_{\M\cap \N}(x_{r_{k_{j}}})<\bar{\mu}_{x^{*}}^{(b_{j}-b)t}d_{\M\cap \N}(x_{a+bt}).
    \]
    Since $\bar{\mu}_{x^{*}} < \mu_{x^{*}}$ we can find a large enough $j$ so that
    \[
    \left(\frac{\bar{\mu}_{x^{*}}}{\mu_{x^{*}}}\right)^{(b_j-b)t}\leq \frac{\mu_{x^{*}}^{a+bt}}{d_{\M\cap \N}(x_{a+bt})}
    \]
    and thus
    \[
    d_{\M\cap \N}(x_{r_{k_{j}}})<\bar{\mu}_{x^{*}}^{(b_{j}-b)t}d_{\M\cap \N}(x_{a+bt})\leq \mu_{x^{*}}^{(b_{j}-b)t}\mu_{x^{*}}^{a+bt} = \mu_{x^{*}}^{r_{k_{j}}}.
    \]
    This contradicts \eqref{eq:counter-ass} so the theorem must hold.
\end{pf}

\begin{rem}
For the case of the method of alternating projections ($\alpha=\alpha_1=\alpha_2=1$), these results coincide with those of \cite{lewis-2008}. In particular, the contraction rate is then given by $\sigma(S_{\T(\bar{x})})=c(\T_{\M(\bar{x})},T_{\N(\bar{x})})$ and the limiting rate is $\gamma(S_{\T(\bar{x})})=c^2(\T_{\M(\bar{x})},T_{\N(\bar{x})})$.
This corresponds to the rates $\cos(\theta_F)$ and $\cos^2(\theta_F)$ where $\theta_F$ is the Friedrichs angle of the corresponding tangent spaces.
\end{rem}

We now show that the faster rate in Theorem \ref{thm:manifold-linear-rate} holds not only in terms of the distance to the intersection,
but also to a point $x^*\in\M\cap\N$.
A similar result can be found in \cite{andersson2013alternating} for the alternating projections method.

\begin{thm}\label{thm:rate-point}
    Let $\M$ and $\N$ be manifolds that satisfy Assumption \ref{ass:regularity} at $\bar{x}\in\M\cap\N$, let the initial
    point $x_{0}\in\reals^n$ be close enough to $\bar{x}$, and let the parameters of the GAP operator $S$ in \eqref{eq:GAP2} satisfy Assumption \ref{ass:alpha} case \ref{ass:alpha-1} or \ref{ass:alpha-2}.
    Further assume that $\M$ and $\N$ satisfy Assumption \ref{ass:regularity} at the limit point $x^*$ of the sequence $(x_k)_{k\in\naturalN}$ generated by the GAP method in Definition~\ref{def:GAP}.
    Then for every $\mu_{x^{*}}\in(\gamma(S_{\T(x^{*})}),1)$, there exists $N\in\naturalN$ such that for all $k\geq N$
    \begin{equation*}
    \|x_{k}-x^*\|\leq\mu_{x^{*}}^{k},\label{eq:limmanifold-point}
    \end{equation*}
    or equivalently
    \begin{equation*}
        \limsup_{k\rightarrow\infty} \|x_{k}-x^*\|^{1/k} \leq \gamma(S_{\T(x^{*})}).
    \end{equation*}
\end{thm}
\begin{pf}
Take any $\mu_{x^{*}}\in(\gamma(S_{\T(x^{*})}),1)$ and let $\bar{\mu} = (\mu_{x^{*}}+\gamma(S_{\T(x^{*})}))/2\le \mu_{x^*}$. Clearly $\bar{\mu}\in(\gamma(S_{\T(x^{*})}),1)$, so we know from Theorem \ref{thm:manifold-linear-rate} that there exists $N$ such that
\begin{equation}\label{eq:copy-linear-rate}
d_{\M\cap\N}(x_k)=\|x_k - \bar{x}_k\|\leq \bar{\mu}^k, \quad \forall k\geq N,
\end{equation}
where $\bar{x}_k\ldef \Proj_{\M\cap\N}(x_k)$.
Pick $c<1$ and $\eta$ so that Theorem~\ref{thm:limitrate1} case \ref{thm:limitrate1:case1}
can be applied for $\bar{x}=x^*$. Since $(x_k)\rightarrow x^*$ there is some $M\geq N$ so that $x_k\in B_{\eta(x^*)}$ for all $k\geq M$ and thus by Theorem~\ref{thm:limitrate1} case \ref{thm:limitrate1:case1}
\begin{equation}\label{eq:copy-slow-rate}
\left\Vert x_{k+1}-\bar{x}_{k}\right\Vert \leq c\left\Vert x_k-\bar{x}_k\right\Vert, \quad \forall k\geq M.
\end{equation}
Using \eqref{eq:copy-linear-rate}, \eqref{eq:copy-slow-rate}, and the triangle inequality we get for $k\geq M$ that
\begin{equation}
\begin{aligned}\
    \|\bar{x}_{k+1}-\bar{x}_{k}\|&\leq\|\bar{x}_{k+1}-x_{k+1}\|+\|{x}_{k+1} -\bar{x}_{k}\| \\
    &\leq \|\bar{x}_{k+1}-x_{k+1}\|+c\|{x}_{k}-\bar{x}_{k}\|\leq \bar{\mu}^{k+1}+c\bar{\mu}^k \\
    &= \bar{\mu}^{k+1}(1+\frac{c}{\bar{\mu}}). 
  \end{aligned}
  \label{eq:bar-bar-inequality}
\end{equation}
By continuity of $\Proj_{\M\cap\N}$ around $x^*$, the point $\bar{x}^*=\lim_{k\rightarrow\infty}\bar{x}_k$ exists.
Using the triangle inequality and \eqref{eq:bar-bar-inequality} we get for $k\geq M$ that
\begin{equation}
\begin{aligned}
    \|\bar{x}_{k}-\bar{x}^*\|&\leq\sum_{i=k}^\infty\|\bar{x}_{i+1}-\bar{x}_{i}\|
    \leq \sum_{i=k}^\infty \bar{\mu}^{i+1}(1+\frac{c}{\bar{\mu}})\\
    & = (1+\frac{c}{\bar{\mu}})\bar{\mu}^{k+1}\sum_{i=0}^\infty \bar{\mu}^{i} \leq (1+\frac{c}{\bar{\mu}})\frac{1}{1-\bar{\mu}}\bar{\mu}^{k+1}
    =\frac{\bar{\mu}+c}{1-\bar{\mu}}\bar{\mu}^{k}.
  \end{aligned}
  \label{eq:bar-star-rate}
\end{equation}
By continuity of $\Proj_{\M\cap\N}$ we also have $x^*=\bar{x}^*$ since $x^*\in\M\cap\N$.
Again, using the triangle inequality, \eqref{eq:copy-linear-rate}, and \eqref{eq:bar-star-rate} we get for $k\geq M$ that
\begin{align*}
    \|x_{k}-x^*\|&\leq \|x_{k}-\bar{x}_k\|+\|\bar{x}_k-x^*\| \leq \bar{\mu}^{k}+\frac{\bar{\mu}+c}{1-\bar{\mu}}\bar{\mu}^{k}=\frac{1+c}{1-\bar{\mu}}\bar{\mu}^{k}.
\end{align*}
Lastly, since $\bar{\mu}<\mu_{x^*}$, there is some $L\geq M$ so that for all $k\geq L$
\begin{equation*}
\|x_{k}-x^*\| \leq \frac{1+c}{1-\bar{\mu}}\bar{\mu}^{k}\leq \mu_{x^*}^k.
\end{equation*}
\end{pf}

We note that the local linear rate $\mu_x^*<\gamma(S_{\T(x^*)})$ is strict,
in the sense that it cannot be improved without adding more assumptions or changing the algorithm.
This follows from the fact that the worst case rate is achieved in the setting of affine sets,
which is covered by this theorem.

To optimize the bound on the convergence rate $\gamma(S_{\T(x^*)})$ in Theorem \ref{thm:rate-point}
in the case where the relative dimensions of the tangent planes are unknown,
it is shown in Theorem \ref{thm:gapiff} that the parameters should be chosen as
\begin{align}
    \alpha=1,\quad\alpha_{1}=\alpha_{2}=\alpha^{*}\ldef \frac{2}{1+\sin{(\theta_{F})}},
\end{align}
where $\theta_{F}$ is the Friedrichs angle between the tangent spaces $\T_{\M(x^*)}$ and $T_{\N(x^*)}$.

\section{Convex sets}

In this section, we show how the convergence results for GAP on manifolds can be extended to GAP on convex sets in some cases. The GAP method is known to converge to some point in the intersection when the sets are closed and convex and the intersection is nonempty, see, e.g., \cite{GAPLS}. The question that remains is the convergence rate.
One way to extend the results in this paper to convex sets is to show that the iterates will eventually behave identically as if the projections were made onto smooth manifolds.
One approach to do this is to partition a convex set into locally smooth manifolds.
This can be done for many convex sets as illustrated in Example~\ref{ex:partitioning}.

\begin{ex}\label{ex:partitioning}\fnote{Example partly stolen from?, add CITE.}
Consider the convex set $C=\{(x,y,z) \mid x^2+y^2\leq z^2, 0\leq z\leq 1\}$. The set can be partitioned into the following five locally smooth manifolds:
$C_1=\inter C,
C_2 = \{(x,y,z) \mid x^2+y^2 = z^2, 0<z<1\},
C_3 = \{(x,y,1) \mid x^2+y^2 < 1\},
C_4 = \{(x,y,1) \mid x^2+y^2 = 1\},
C_5 = \{(0,0,0)\}.$
\end{ex}

There is plenty of literature on this type of identification of surfaces. For example, in \cite{act_ident_15} the authors study the Douglas--Rachford algorithm for partially smooth functions. However, the assumptions do not generally apply to convex feasibility problems since all reformulations into the framework will either be non-smooth or have vanishing gradients at the boundaries.

For the case of alternating projections on convex sets, the projections will always lie on the boundary of the sets until the problem is solved. The local convergence rate therefore follows trivially if the boundaries of these sets satisfy the required regularity assumptions at the intersection.
However, this is not the case for GAP in general because of the (over)-relaxed projections.
Even in cases of polyhedral sets,
identification of affine sets is not guaranteed as we show with an example in Section \ref{sec:counter-example}.
We therefore present our results under some smoothness assumptions and for a slightly restricted set of parameters. This set of parameters does include the set of parameters that optimize the rate in Theorem \ref{thm:rate-point}.



\begin{lem}\label{lem:ball}%
    Let $A$ be a closed solid convex set in $\realR^n$ with a $\C^2$ smooth boundary around $\bar{x}\in\bd A$.
    Then there exists a $\delta>0$ such that for all $x\in\Ball_\delta(\bar{x})\setminus A$
    \begin{align*}
        \Proj_A^\alpha x\in \inter A,\,\, \forall \alpha \in (1,2].
    \end{align*}
\end{lem}
\begin{pf}
    As noted in Remark~\ref{rem:smooth-normal}, smoothness of $\bd A$ at $\bar{x}\in\bd A$ implies that there exists a neighborhood $U$ of $\bar{x}$ for which the outwards facing normal vector $n(x)$ with $\|n(x)\|=1$ is unique for all $x\in\bd A\cap U$ and that the normal $n(x)$ is continuous around $\bar{x}$.
    Since $A$ is solid and $\bd A$ is smooth at $\bar{x}$, there is some $\zeta>0$ so that $\bar{x}-\beta  n(\bar{x})\in \inter A$ for all $\beta\in(0,\zeta]$.
    We assume without loss of generality that $\zeta<1$.
    We can now create an open ball with radius $\delta$ such that
    \begin{equation}\label{eq:ball-interior}
        \Ballo_\delta(\bar{x}-\beta  n(\bar{x}))\subset \inter A.
    \end{equation}
    From continuity of $n(x)$ we have that there exists $\epsilon'>0$ such that for all $x\in \bd A$
    \begin{equation}\label{eq:ball-open-epsilon}
        \|x-\bar{x}\|\leq \epsilon' \Rightarrow \|n(x)-n(\bar{x})\|\leq \delta.
    \end{equation}
    Now pick $0<\epsilon<\min(\delta(1-\beta ), \beta , \epsilon')$. By the triangle inequality, for all $x \in\Ball_\epsilon(\bar{x})\cap \bd A$,
    \begin{align*}
        \|(x-\beta n(x))-(\bar{x}-\beta n(\bar{x}))\| &\leq \|x-\bar{x}\|+\beta\|n(x)-n(\bar{x}))\| \leq \epsilon + \beta \delta  < \delta(1-\beta )+\beta \delta = \delta.
    \end{align*}
    Using this and \eqref{eq:ball-interior},
    \begin{equation}\label{eq:ball-interior-epsilon}
        x-\beta n(x) \in \inter A\,, \forall x \in \Ball_\epsilon(\bar{x})\cap \bd A.
    \end{equation}
    Moreover, by convexity of $A$ and non-expansiveness~\cite[Prop. 4.16]{bauschkeCVXanal} of the projection
    \begin{equation}\label{eq:ball-projection-close}
        \Proj_A(x)\in\Ball_\epsilon(\bar{x}), \forall x\in\Ball_\epsilon(\bar{x}).
    \end{equation}
    Hence, by \eqref{eq:ball-interior-epsilon}, \eqref{eq:ball-projection-close}, and since $\Proj_A(x)\in\bd(A)$ for $x\not\in A$ we have
    \begin{equation}
    \Proj_A(x) - \beta n(\Proj_A(x))\in \inter A,\, \forall x\in \Ball_\epsilon(\bar{x}) \setminus A.
    \end{equation}
    Moreover, the projection operator satisfies
    \begin{align*}
        n(\Proj_A(x)) = \frac{x-\Proj_A(x)}{\|x-\Proj_A(x)\|},
    \end{align*}
    for $x\not\in A$~\cite[Prop. 6.47]{bauschkeCVXanal}.
    By the definition of the relaxed projection
    we therefore have for $x\in\Ball_\epsilon(\bar{x}) \setminus A$ that $\Proj_A^\alpha(x) = \Proj_A(x) - (\alpha-1)\|\Proj_A(x)-x\|n(\Proj_A(x))$.
    Noting that since $\alpha\in(1,2]$ we have
    \[
    0<(\alpha-1)\|\Proj_A(x)-x\|\leq \epsilon < \beta < 1,
    \]
    which implies that $\Proj_A^\alpha(x)$ is a strict convex combination
    between
    $\Proj_A(x)\in A$ and $\Proj_A(x) - \beta n(\Proj_A(x))\in \inter A$, i.e.,
    \[
        \Proj_A^\alpha(x) = \gamma \Proj_A(x) + (1-\gamma) (\Proj_A(x)-\beta n(\Proj_A(x))),
    \]
    where $\gamma\ldef1-(\alpha-1)\|\Proj_A(x)-x\|/\beta\in(0,1)$
    and therefore $\Proj_A^\alpha(x)\in \inter A$.
\end{pf}

\subsection{Examples}
In this section, we present some results on when the rate in Theorem \ref{thm:rate-point} can be applied to convex sets. We say for a convex set $A\in\reals^n$ that the algorithm has \emph{identified} a manifold $\M\subset A$ at some iteration $k$ if subsequent iterations would be identical when the set $A$ is replaced with $\M$. We partition a smooth convex set $A$ into two parts $\bd A$ and $\inter A$, and show that either $\bd A$ or $\inter A$ is identified.
%
\begin{ass}[regularity of convex sets at solution]\label{ass:regularity-convex}%
    Let $A,B\subseteq\reals^n$ be two closed convex sets with $x^*\in A \cap B$.
    Assume that at least one of the following holds
    \begin{enumerate}[label=C\arabic*.,ref=C\arabic*]
        \item $x^*\in\bd A\cap\bd B$ and $(\bd A,\bd B)$ satisfies Assumption
        \ref{ass:regularity} at the point $x^*$, \label{ass:regularity-convex-bd}
        \item $x^*\in \inter A\cap \bd B$ where $\bd B$ is $\mathcal{C}^2$-smooth around $x^*$, \label{ass:regularity-convex-inter-1}
        \item $x^*\in \bd A\cap \inter B$ where $\bd A$ is $\mathcal{C}^2$-smooth around $x^*$, \label{ass:regularity-convex-inter-2}
        \item $x^*\in \inter A \cap \inter B$.\label{ass:regularity-convex-inter-3}
    \end{enumerate}
\end{ass}

We now introduce a definition of $S_{\T(x^*)}$ in the setting of convex sets to simplify the following statements on convergence rates.
\begin{defin}
    For the feasibility problem $(A,B)$ involving two closed convex sets $A,B\subseteq\reals^n$ that satisfy Assumption \ref{ass:regularity-convex} at a point $x^*\in A\cap B$, we define
    \begin{equation*}
        S_{\T(x^*)} \ldef (1-\alpha)I+\alpha\Proj_{\T_{\M}(x^*)}^{\alpha_{2}}\Proj_{\T_{\N}(x^*)}^{\alpha_{1}}
    \end{equation*}
    where we let
    \begin{equation*}
        \M \ldef \begin{cases}
        \bd A & \textrm{ if } x^*\in\bd A\\
        \inter A & \textrm{ if } x^*\in\inter A
        \end{cases},\quad
        \N \ldef \begin{cases}
        \bd B & \textrm{ if } x^*\in\bd B\\
        \inter B & \textrm{ if } x^*\in\inter B.
        \end{cases}
    \end{equation*}
\end{defin} 

If $x^*\in \inter A$ in this definition,
we get tangent space $\T_{\M}(x^*)=\realR^n$ and projection operator
$\Proj_{\T_{\M}(\bar{x})}^{\alpha_{2}}=I$. The same holds for $x^*\in \inter B$.
The corresponding rate $\gamma(S_{\T(x^*)})$ then reduces to one of $(1-\alpha_2)$, $(1-\alpha_1)$ or $(1-\alpha_1)(1-\alpha_2)$ according to Theorem \ref{thm:gap-all-eigs}.

\begin{thm}\label{thm:convex-rate}
    Consider $(A,B)$, let $A,B\subseteq\reals^n$ be solid closed convex sets with $A\cap B\neq \emptyset$, and let $\alpha=1, 1<\alpha_1,\alpha_2<2$ in the GAP algorithm in Definition~\ref{def:GAP}. 
    Then the iterations converge to some point $x^*\in A\cap B$.
    If the sets $A$ and $B$ satisfy Assumption \ref{ass:regularity-convex} at the point $x^*$, then
    either the problem is solved in finite time,
    or eventually the algorithm will identify the sets $\bd A$ and $\bd B$ and converge R-linearly with any rate $\mu\in (\gamma(S_{\T(x^{*})}),1)$ to $x^*\in\bd A\cap \bd B$.
\end{thm}
\begin{pf}
    We know that $x_{k}\rightarrow x^*$ for some point $x^*$ due to closed convexity of $A$ and $B$~\cite[Prop. 3]{GAPLS}.
    We first show that the problem is solved in a finite number of iterations unless $x^*\in\bd A\cap \bd B$.

    Assume $x^*\in \inter A\cap  \inter B$. Then there is some open ball around $x^*$ that is contained in $A\cap B$.
    By convergence of $(x_k)_{k\in\naturalN}$, there is some $k$ such that $x_{k}$ is in this ball, and we have convergence in finite time.

    Assume $x^*\in \bd A\cap \inter B$. Let $\delta$ be such that Lemma \ref{lem:ball} is satisfied for $(A,x^*$)
    and so that $\Ball_\delta(x^*)\subset B$.
    Then there is a $k\in\naturalN$ such that $x_k\in \Ball_{\delta}(x^*)$.
    If $x_{k}\in A\cap B$ the problem is solved in finite time.
    If not, then $x_{k}\in B\setminus A$, so trivially $\Proj_B^{\alpha_1}x_k=x_k$, and by Lemma \ref{lem:ball} we get $x_{x+1}=\Proj_A^{\alpha_2}x_k\in\inter A$.
    By non-expansiveness of $\Proj_A^{\alpha_2}\Proj_B^{\alpha_1}$,
    we have $x_{k+1}\in\Ball_{\delta}(x^*)\subset B$,
    so $x_{k+1}\in A\cap B$, and the problem is solved in finite time.

    Assume $x^*\in \inter A\cap \bd B$ and let $\delta$ be such that Lemma \ref{lem:ball} is satisfied for $(B,x^*$), and so that $\Ball_\delta(x^*)\subset A$.
    Eventually $x_k\in \Ball_{\delta}(x^*)$ for some $k\in\naturalN$.
    If $x_k\in B$ the problem is solved.
    If not, then $x_k\in A\setminus B$,
    but then $\Proj^{\alpha_1}_Bx_k\in B$ by Lemma~\ref{lem:ball}.
    Again, by non-expansiveness of $\Proj_B^{\alpha_1}$ we have $\Proj^{\alpha_1}_Bx_k\in\Ball_{\delta}(x^*)\subset A$ so $x_{k+1}=\Proj^{\alpha_1}_Bx_k\in A\cap B$ and the problem is solved in finite time.

    Now consider the case where $x^*\in \bd A \cap \bd B$. Choose $\delta_A$ and $\delta_B$ so that Lemma \ref{lem:ball} is satisfied for $(A,x^*)$ and $(B,x^*)$ respectively and let $\delta=\min(\delta_A,\delta_B)$.
    Since $x_k\rightarrow x^*$ there exists $N\in\naturalN$ such that $x_k\in\Ball_\delta(x^*)$ for all integer $k>N$. By Lemma \ref{lem:ball},
    we then have $x_{k+1}\in A$. If $x_{k+1}\in A\cap B$, the problem is solved in finite time. If not, $x_{k+1}\in A\setminus B$.
    Now consider any integer $j>N$ such that $x_{j}\in A\setminus B$ with $x_j\in\Ball_\delta(x^*)$. The first projection $\Proj_B^{\alpha_1}(x_{j})$ is equivalent to projecting onto the manifold $\bd B$, and by Lemma \ref{lem:ball}, we have $\Proj_B^{\alpha_1}(x_{j})\in B$.
    Either $\Proj_B^{\alpha_1}(x_{j})$ is also in $A$, in which case the problem is solved in finite time, or the second projection $\Proj_A^{\alpha_2}\Proj_B^{\alpha_1}(x_{j})$ is equivalent to projecting onto the manifold $\bd A$. By Lemma \ref{lem:ball}, we get $x_{j+1}\in A$.
    Therefore, we either we have $x_{j+1}\in A\cap B$, in which case we have a solution in finite time, or we have $x_{j+1}\in A\setminus B$.
    By recursion over $j>N$, we see that either the problem is solved in finite time,
    or $x_{j+1}\in A\setminus B$ for all $j>N$,
    in which case each projection onto the sets is equivalent to projecting onto their boundaries,
    i.e., the algorithm has identified the manifolds. The rate then follows directly from Theorem \ref{thm:rate-point}.
\end{pf}
%

\begin{thm}\label{thm:rate-convex-sets}
    Let $A\subseteq\reals^n$ be a solid closed convex set and $B\subseteq\reals^n$ be an affine set such that $A\cap B\neq \emptyset$.
    Then $x_k\rightarrow x^*$ for some point $x^*\in A\cap B$ for the GAP algorithm in Definition~\ref{def:GAP}. 
    If the sets $A$ and $B$ satisfy Assumption \ref{ass:regularity-convex} at $x^*$, then
    the GAP method converges R-linearly with any rate $\mu\in(\gamma(S_{\T(x^{*})}),1)$ to $x^*$.
\end{thm}
\begin{pf} This proof is similar to that of Theorem \ref{thm:convex-rate}.
    The sequence $(x_k)_{k\in\naturalN}$ converges to some $x^*\in A\cap B$ by convexity of the sets.
    First assume that $x^*\in\inter A$.
    Then, since $x_k\rightarrow x^*$ there exists $N$ such that $x_j\in A$ for all $j>N$.
    The problem is then locally equivalent to that of $(\realR^n,B)$, i.e., two subspaces.

    If $x^*\in\bd A$, then let $\delta$ be such that Lemma \ref{lem:ball} is satisfied for $(A,x^*)$. Then by convergence to $x^*$, eventually $x_j\in\Ball_\delta(x^*)$ for all $j>N$.
    If $\Proj_B^{\alpha_1}x_j\not\in A$ then $x_{j+1}\in \inter A$ by Lemma \ref{lem:ball}. If $\Proj_B^{\alpha_1}x_j\in A$, then $x_{j+1}\in A$ by the definition of projection.
    So $x_{j+1}\in A$ for all $j>N$.

    If also $\Proj_B^{\alpha_1}x_l\in A$ for some $l>j>N$,
    then since both $x_l$ and $x_{l-1}$ are in $A$,
    we have $x_l-x_{l-1}\in \Normal_B(\Proj_Bx_{l-1})$.
    From convexity of $A$ we know that the straight line segment between $x_l$ and $x_{l-1}$ must be contained in $A$,
    so all subsequent iterations must be on this line segment.
    But then $\Proj_Bx_l=x^*$ and by assumption $x^*\in \bd A$, so convexity of $A$ implies that the whole segment must be in $\bd A$.
    The algorithm has thus identified $\bd A$ and $B$.

    Otherwise, $\Proj_B^{\alpha_1}x_j \not\in A$ for all $j>k$,
    and the projection $\Proj_A^{\alpha_2}(\Proj_B^{\alpha_1})x_j$ is equivalent to projecting onto $\bd A$, i.e, the algorithm has identified $\bd A$ and $B$. The rate then follows from Theorem \ref{thm:rate-point} since $B$ is a smooth manifold.
\end{pf}

Next, we introduce some regularity properties of convex sets and show how they relate to the regularity of the manifolds corresponding to their boundaries.

\begin{defin}[subtransversality of sets]\cite[Thm. 1 (ii)]{kruger2018set}\label{def:sr}\\
    Two sets $C\subseteq\reals^n$ and $D\subseteq\reals^n$ are \emph{subtransversal} at $x^*\in C\cap D$ if there exist $\alpha>0$ and $\delta>0$ such that
    \begin{align}\label{eq:sr}
        \alpha \dist_{C\cap D}(x)\leq\max\{\dist_{C}(x), \dist_{D}(x)\}\quad  \forall x\in\Ball_\delta(x^*).
    \end{align}
     $\sreg[C,D](x^*)$ is defined as the exact upper bound of all $\alpha$ such that \eqref{eq:sr} holds.
\end{defin}

\begin{defin}[transversality of sets]\cite[Thm. 1 (ii)]{kruger2018set}\\
    Two sets $C\subseteq\reals^n$ and $D\subseteq\reals^n$ are \emph{transversal} at $x^*\in C\cap D$ if there exists $\alpha>0$ and $\delta>0$ such that
    \begin{align}\label{eq:reg}
        \alpha \dist_{(C-x_1)\cap (D-x_2)}(x)\leq&\max\{\dist_{C-x_1}(x), \dist_{D-x_2}(x)\}\quad {\hbox{ for all }} x\in\Ball_\delta(x^*), x_1,x_2\in\Ball_\delta(0).
    \end{align}
     $\reg[C,D](x^*)$ is defined as the exact upper bound of all $\alpha$ such that \eqref{eq:reg} holds.
     Equivalently, the sets $C$ and $D$ are transversal at $x^*$ if $\Normal_C(x^*)\cap(-\Normal_D(x^*))=\{0\}$~\cite[Thm. 2 (v)]{kruger2018set}.
\end{defin}
The transversality condition $\Normal_C(x^*)\cap(-\Normal_D(x^*))=\{0\}$ for two sets $C$ and $D$
coincides with Definition \ref{def:substransversal-manifolds} of transversality when the sets are smooth manifolds,
since the normal cones are linear subspaces in this case~\cite{paul1947finite}.

\begin{defin}[acute and obtuse intersection]
    For two solid, closed, convex sets $A,B\subseteq\reals^n$ with smooth boundaries, we say that the intersection is \emph{acute} at a point $x^*\in\bd A\cap \bd B$ if $\langle v_1, v_2\rangle \leq 0$, where $v_1,v_2$ are the unique vectors such that $v_1\in\Normal_A(x^*), v_2\in\Normal_B(x^*), \|v_1\|=\|v_2\|=1$. Conversely, we say that the intersection is \emph{obtuse} if $\langle v_1, v_2\rangle > 0$.
\end{defin}
Note that \emph{acute} and \emph{obtuse} refer to the shape of the intersection, and not the angle between the normals, for which the property is reversed.
\begin{lem}\label{lem:set-transversality}
Let $A,B$ be solid, closed and convex sets in $\realR^n$ with boundaries $\bd A,\bd B$
that satisfy Assumption~\ref{ass:regularity} at some point $x^*\in\bd A,\bd B$
and assume that $\T_{\bd A}(x^*)\neq \T_{\bd B}(x^*)$.
Let $\theta_F\in(0,\pi/2]$ be defined via $\cos(\theta_F)=c(\bd A, \bd B, x^*)$.
Then
\begin{enumerate}[label=\arabic*.,ref=\arabic*]
    \item\label{lem:set-transversality:casebdAbdB}
        the manifolds $\bd A$ and $\bd B$ are transversal at $x^*$,
    \item\label{lem:set-transversality:caseAB}
        the sets $A$ and $B$ are transversal at $x^*$, i.e. $\Normal_{A}(x^*)\cap(-\Normal_{B}(x^*))=\{0\}$,
    \item\label{lem:set-transversality:AB}
    the sets $A$ and $B$ are subtransversal at $x^*$ and
    the following inequalities hold
    \begin{align*}
        \reg[A,B](x^*)\leq\sreg[A,B](x^*) \leq \begin{cases}
        \sin(\theta_F/2) \quad\textrm{ if } (A,B) \textrm{ acute at } x^*\\
        \cos(\theta_F/2) \quad\textrm{ if } (A,B) \textrm{ obtuse at } x^*,
        \end{cases}
    \end{align*}
    \item\label{lem:set-transversality:bdAB}
    $\sin(\theta_F/2)=\reg[\bd A, \bd B](x^*)$.
    Furthermore, if the intersection of $(A,B)$ is acute at $x^*$ then
    \begin{align*}
      \sin(\theta_F/2) =\reg[\bd A, \bd B](x^*)=\reg[A,B](x^*)=\sreg[A,B](x^*)
    \end{align*}
    otherwise
    \begin{align*}
        \cos(\theta_F/2)=\reg[A,B](x^*)=\sreg[A,B](x^*).
    \end{align*}
\end{enumerate}
\end{lem}
\begin{pf}
    The proofs follow the definitions and results on (sub-)transversality of general sets from~\cite{kruger2006regularity}.

    \ref{lem:set-transversality:casebdAbdB}: From smoothness of the manifolds $\bd A, \bd B$, the corresponding normals are lines and trivially $\Normal_{\bd B}(x^*)=-\Normal_{\bd B}(x^*)$.
    Moreover, since $\T_{\bd A}(x^*)\neq \T_{\bd B}(x^*)$ we have $\Normal_{\bd A}(x^*)\neq\Normal_{\bd B}(x^*)$, and therefore
    $\Normal_{\bd A}(x^*)\cap(-\Normal_{\bd B}(x^*))=\{0\}.$

    \ref{lem:set-transversality:caseAB}: The normals to the sets $A,B$ at a point on their boundaries $x^*$ satisfy $\Normal_{\bd A}(x^*)=\Normal_A(x^*)\cup(-\Normal_A(x^*))$ and correspondingly for $B$. Hence, $\Normal_A(x^*)\subset\Normal_{\bd A}(x^*)$ and $-\Normal_B(x^*)\subset\Normal_{\bd B}(x^*)$,
    so from case~\ref{lem:set-transversality:casebdAbdB} it follows that
    $\Normal_{A}(x^*)\cap(-\Normal_{B}(x^*))=\{0\}$.

    \ref{lem:set-transversality:AB}:
    The first inequality follows directly from~\cite[Thm. 4 (i)]{kruger2018set}.
    For the second inequality, let $v_1\in\Normal_A(x^*)$,$v_2\in\Normal_B(x^*)$ be the unique vectors with
    $\|v_1\|=\|v_2\|=1$, and define
    $w=(v_1+v_2)/\|v_1+v_2\|$.
    From case~\ref{lem:set-transversality:caseAB}, we see that $v_1\neq -v_2$ and thus $\langle v_1,v_2\rangle>-1$.
    Thus $\langle w, v_1\rangle = (\langle v_1, v_2\rangle + 1)/\|v_1+v_2\| > 0$ and similarly $\langle w, v_2\rangle>0$.
    Since $A,B$ are convex sets, $\T_{A}(x^*)+\{x^*\}$ and $\T_{B}(x^*)+\{x^*\}$ are separating hyperplanes to the corresponding sets,
    and it follows from $\langle w, v_1\rangle>0, \langle w, v_2\rangle>0$ that $x^*+\beta w$ is separated from the sets $A$ and $B$ when $\beta>0$, i.e. $x^*+\beta w\not\in A\cup B$ for $\beta>0$. Moreover, by definition of $w$, we have $w\in N_A(x^*)+N_B(x^*)\subset N_{A\cap B}(x^*)$ where the second inclusion holds trivially for convex sets. We can therefore conclude that $\Proj_{A\cap B}(x^*+\beta w)=x^*$,
    and therefore
    \begin{align}\label{eq:distAB}
        \dist_{A\cap B}(x^*+\beta w)=\beta\|w\|=\beta.
    \end{align}
    We now calculate an expression for $\dist_{A}(x^*+\beta w)$.
    Since $x^*+\beta w\not \in A$, the projection onto $A$ is locally equivalent to projecting onto the smooth manifold $\bd A$.
    From Lemma~\ref{lem:proj-manifold} we get with series expansion around $x^*$ that
    \begin{align*}
        \Proj_{\bd A}(x^*+\beta w) = \Proj_{\bd A}(x^*) + \Proj_{\T_{\bd A}(x^*)}(\beta w) + O(\beta^2),
    \end{align*}
    where $\Proj_{\bd A}(x^*)=x^*$.
    The projection of $w=(v_1+v_2)/\|v_1+v_2\|$ onto $\T_{\bd A}(x^*)$ is given by
    \begin{align*}
        \Proj_{\T_{\bd A}(x^*)}(w) &= w - \frac{\langle v_1, w\rangle}{\|v_1\|^2}v_1= w - \langle v_1, w\rangle v_1
    \end{align*}
    and the distance $d_{A}(x^*+\beta w)$ is therefore
\begin{equation}
    \begin{aligned}
        d_{A}(x^*+\beta w) &= \|\Proj_{\bd A}(x^*+\beta w)-(x^*+\beta w)\|= \|\beta\Proj_{\T_{\bd A}(x^*)}(w)-\beta w+O(\beta^2)\|\\
        &  = \|\beta\langle v_1, w\rangle v_1 - O(\beta^2)\|  = \beta\|\frac{1+\langle v_1, v_2\rangle}{\|v_1+v_2\|}v_1 - O(\beta)\|,
      \end{aligned}
      \label{eq:distA}
    \end{equation}
    and in the same way for $B$: $d_{B}(x^*+\beta w)=\beta\|\frac{1+\langle v_1, v_2\rangle}{\|v_1+v_2\|}v_2 - O(\beta)\|$.

    By Definition \ref{defin:friedrichs} that defines the Friedrichs angle and Definition~\ref{defin:cmn}, we conclude that
    \begin{align*}
        \cos(\theta_F)&=c(\bd A,\bd B, x^*)=c(\T_{\bd A}(x^*), \T_{\bd B}(x^*))=c((\T_{\bd A}(x^*))^\perp, (\T_{\bd B}(x^*))^\perp),
    \end{align*}
    where the last equality can be found, e.g., in \cite[Def. 3]{kruger2018set}.
    Since $(\T_{\bd A}(x^*))^\perp=\Normal_A(x^*)\cup (-\Normal_A(x^*))=\{\beta v_1 \mid \beta\in\realR\}$, and similarly for $B$, Definition \ref{defin:friedrichs} results in that $\cos(\theta_F)=\max\{\langle v_1,v_2\rangle, -\langle v_1,v_2\rangle\}$, i.e.
    \begin{align*}
        \langle v_1, v_2 \rangle =
        \begin{cases}
         -\cos(\theta_F) & \textrm{ if } \langle v_1, v_2 \rangle \leq 0\\
         \cos(\theta_F) & \textrm{ if } \langle v_1, v_2 \rangle \geq 0.
        \end{cases}
    \end{align*}
    Thus by definition of $\sreg[A,B](x^*)$, \eqref{eq:distAB} and \eqref{eq:distA}
    \begin{align*}
        \sreg[A,B](x^*) &\leq \lim_{\beta\rightarrow 0^+}\frac{\max(\dist_A(x^*+\beta w),\dist_B(x^*+\beta w))}{\dist_{A\cap B}(x^*+\beta w)}\\
        &=\lim_{\beta\rightarrow 0^+}\max_{i\in\{1,2\}}\|\frac{1+\langle v_1, v_2\rangle}{\|v_1+v_2\|}v_i - O(\beta)\|\\
        &=\frac{1+\langle v_1, v_2\rangle}{\sqrt{\|v_1\|^2+2\langle v_1, v_2 \rangle + \|v_2\|^2}}\\
        &=\begin{cases}
        \frac{1-\cos(\theta_F)}{\sqrt{2-2\cos(\theta_F)}}=\sqrt{1-\cos(\theta_F)}/\sqrt{2}=\sin(\theta_F/2) \,\,\,\textrm{ if } \langle v_1, v_2 \rangle \leq 0\\
        \frac{1+\cos(\theta_F)}{\sqrt{2+2\cos(\theta_F)}}=\sqrt{1+\cos(\theta_F)}/\sqrt{2}=\cos(\theta_F/2) \,\,\,\textrm{ if } \langle v_1, v_2 \rangle \geq 0.
        \end{cases}
    \end{align*}
    \ref{lem:set-transversality:bdAB}:
    By~\cite[Prop. 8]{kruger2018set}
    \begin{align*}
        \reg_\textrm{a}[C, D](x) = \sup_{\substack{n_1\in\Normal_C(x),\, n_2\in\Normal_D(x) \\ \|n_1\|=\|n_2\|=1}} -\langle n_1, n_2\rangle,
    \end{align*}
    where $\reg_\textrm{a}[C, D](x)$ satisfies $\reg_\textrm{a}[C, D](x^*)+2(\reg[C, D](x^*))^2=1$.

    Since $\bd A,\bd B$ are smooth manifolds, this results in $\reg_\textrm{a}[\bd A, \bd B](x^*)=\cos(\theta_F)$
    by Definition~\ref{defin:friedrichs},
    since $\Normal_{\bd A}(x^*)=-\Normal_{\bd A}(x^*)$ and equivalently for $\bd B$.
    Thus, since $\theta_F\in[0,\pi/2]$ and $\reg[\bd A, \bd B](x^*)\geq0$ holds by definition, we have
    $\reg[\bd A, \bd B](x^*)=\sqrt{(1-\cos(\theta_F))/2}=\sin(\theta_F/2)$ for all $\theta_F\in[0,\pi/2]$.

    For $\reg[A,B](x^*)$ we use the same result, but the unit normal vectors are unique in this case.
    When $\langle v_1, v_2\rangle\leq0$ we have $\langle v_1, v_2\rangle = -\cos(\theta_F)$ by definition of $\theta_F$.
    We therefore get $\reg_\textrm{a}[A, B]=\cos(\theta_F)$ and thus $\reg[A,B](x^*)=\sqrt{(1-\cos(\theta_F))/2}=\sin(\theta_F/2)$.

    In the same way, when $\langle v_1, v_2\rangle\geq0$ we have $\langle v_1, v_2\rangle = \cos(\theta_F)$,
    so $\reg_\textrm{a}[A, B]=-\cos(\theta_F)$ and $\reg[A,B](x^*)=\sqrt{(1+\cos(\theta_F))/2}=\cos(\theta_F/2)$.

    But we always have $\reg[A,B]\leq \sreg[A,B]$~\cite[Thm. 4 (i)]{kruger2018set}, so together with case \ref{lem:set-transversality:AB} we see that $\sreg[A,B](x^*)$ is bounded both above and below by
    \begin{align*}
        \sin(\theta_F/2)& \quad\textrm{ if } \langle v_1, v_2 \rangle \leq 0\\
        \cos(\theta_F/2)& \quad\textrm{ if } \langle v_1, v_2 \rangle \geq 0,
    \end{align*}
    which concludes the proof.
\end{pf}
\begin{rem}
    The regularity constants in Lemma~\ref{lem:set-transversality} case \ref{lem:set-transversality:bdAB} are continuous with respect to the normals as they approach the limit between acute and obtuse since $\langle v_1,v_2\rangle \rightarrow 0 \Rightarrow \theta_F\rightarrow \pi/2$ and $\sin(\pi/4)=\cos(\pi/4)=1/\sqrt{2}$.
\end{rem}
The rates presented so far are stated either as a property of the operator $S_{\T(x^*)}$ or as a function of the Friedrichs angle $\theta_F$ between tangent planes at the intersection.
In previous work on alternating projections and similar algorithms for convex and non-convex sets,
the rates are often stated as a function of a linear regularity constant~\cite{lewis2009local,bauschke2014MARP}.
We now state the rate found by choosing the optimal relaxation parameters~\eqref{eq:optpar} in terms of linear regularity.
\begin{thm}\label{thm:regularity}
    Let $A,B$ be two solid, closed, and convex sets in $\realR^n$. Let $x^*\in A\cap B$ be the limit point of the sequence $(x_k)_{k\in\naturalN}$ generated by the GAP algorithm in Definition~\ref{def:GAP}, and assume that
    \begin{enumerate}
        \item $x^*\in\bd A\cap\bd B$
        \item $\bd A$ and $\bd B$ satisfy Assumption \ref{ass:regularity} at the point $x^*$.
    \end{enumerate}
    Then the sets are $\hat{\kappa}$-linearly regular, i.e., there exist $\delta>0$ and $\hat{\kappa}>0$ such that
    \begin{equation}\label{eq:regularity-kappa}
        \dist_{A\cap B}(x) \leq \hat{\kappa} \max(\dist_A(x),\dist_B(x)), \quad \forall x\in\Ball_\delta(x^*).
    \end{equation}
    Let $\kappa$ be the lower limit of all such $\hat{\kappa}$ and assume that $\kappa\geq\sqrt{2}$, then the GAP algorithm with parameters
    \begin{equation}
        \alpha=1, \quad\alpha_1=\alpha_2=2\left(\frac{\kappa}{\sqrt{\kappa^2-1}+1}\right)^2
    \end{equation}
    will converge to $x^*$ with R-linear rate $\mu$ for any $\mu\in(\gamma,1)$, where
    \begin{equation}\label{eq:kappa-rate}
        \gamma=\left(\frac{\sqrt{\kappa^2-1}-1}{\sqrt{\kappa^2-1}+1}\right)^2=1-4\frac{\sqrt{\kappa^2-1}}{\kappa^2+2\sqrt{\kappa^2-1}}.
    \end{equation}
\end{thm}
\begin{pf}
    Existence of a limit point $x^*$ for convex sets follows from the previous results or~\cite{GAPLS}.
    First assume that $T_{\bd A}(x^*)=T_{\bd B}(x^*)$.
    Then Assumption \ref{ass:reg-intersection} along with a dimensionality argument imply that $\bd A=\bd B$ in some neighborhood of $x^*$. It must therefore be that either $A\cap B=A=B$ or $A\cap B = \bd A \cap \bd B$ in some neighborhood of $x^*$.
    The problem is then trivial, but $\dist_{A\cap B}(x)=\dist_A(x)=\dist_B(x)$ for all $x\in\Ball_\delta(x^*)$, so $\kappa=1$. This trivial case is not covered by the result as we assume $\kappa\geq\sqrt{2}$.

    Now assume instead that $T_{\bd A}(x^*)\neq T_{\bd B}(x^*)$. The sets ($A,B$) are therefore transversal by Lemma \ref{lem:set-transversality} case \ref{lem:set-transversality:caseAB}, and since $\Normal_A(x^*)\neq\Normal_B(x^*)$, we have $\theta_F>0$.
    Since $1/\kappa=\sreg[A,B]\leq1/\sqrt{2}$ we have by Lemma~\ref{lem:set-transversality} case~\ref{lem:set-transversality:bdAB} that
    \begin{align*}
        1/\kappa=\reg[\bd A, \bd B]=\sreg[A,B]=\sin(\theta_F/2).
    \end{align*}
    The optimal parameters \eqref{eq:optpar} are therefore, with $\theta_F=2\arcsin(1/\kappa)$
    \begin{align*}
      \alpha_1=\alpha_2=\frac{2}{1+\sin(\theta_F)} = \frac{2}{1+\sin(2\arcsin(1/\kappa))} =
        2\left(\frac{\kappa}{\sqrt{\kappa^2-1}+1}\right)^2\in[1,2).
    \end{align*}
    By Theorem \ref{thm:rate-convex-sets} and Theorem \ref{thm:gapiff}, the convergence to $x^*$ is R-linear with rate $\mu$ for any $\mu\in(\gamma(S_{\T(x^{*})}),1)$ where
    \begin{align*}
        \gamma(S_{\T(x^{*})}),1)&=\frac{1-\sin(\theta_F)}{1+\sin(\theta_F)}=\frac{1-\sin(2\arcsin(1/\kappa))}{1+\sin(2\arcsin(1/\kappa))}=
        \left(\frac{\sqrt{\kappa^2-1}-1}{\sqrt{\kappa^2-1}+1}\right)^2=1-4\frac{\sqrt{\kappa^2-1}}{\kappa^2+2\sqrt{\kappa^2-1}}.
    \end{align*}
\end{pf}

\begin{rem}
    The regularity parameter satisfies $\kappa\in[1,\infty]$ under the assumptions of Theorem~\ref{thm:regularity}.
    In particular, for ill-conditioned problems, i.e., large $\kappa$, the rate above approaches $\gamma\approx 1-\frac{4}{\kappa}$.
    This can be compared to the rate of alternating projections of $\gamma=1-\frac{4}{\kappa^2}$ as found in \cite{lewis2009local} under linear regularity assumptions for non-convex sets, which is worse.
    The contraction rate for the Douglas--Rachford algorithm, presented in \cite{Luke2020} for general convex sets is $\sqrt{1-\kappa^{-2}}$, which can be approximated for large $\kappa$ by $1-\frac{1}{2\kappa^2}$.

\end{rem}
\begin{thm}\label{thm:regularity-global}
    Let $A,B$ be two solid, closed, and convex sets in $\realR^n$ that satisfy Assumption~\ref{ass:regularity-convex} at every point $x^*\in A \cap B$.
    Assume that there is a $\hat{\kappa}>0$ such that the sets $A,B$ are $\hat{\kappa}$-linearly regular at every point $x^*\in A\cap B$, i.e.,
    for every $x^*$ there exists $\delta_{x^*}>0$ such that
    \begin{equation}\label{eq:regularity-kappa2}
        \dist_{A\cap B}(x) \leq \hat{\kappa} \max(\dist_A(x),\dist_B(x)), \quad \forall x\in\Ball_{\delta_{x^*}}(x^*).
    \end{equation}
    Let $\kappa=\max(\hat{\kappa}, \sqrt{2})$, then the GAP algorithm with parameters
    \begin{equation}
        \alpha=1, \quad\alpha_1=\alpha_2=2\left(\frac{\kappa}{\sqrt{\kappa^2-1}+1}\right)^2
    \end{equation}
    will converge to $x^*$ with R-linear rate $\mu$ for any $\mu\in(\gamma,1)$, where
    \begin{equation}\label{eq:kappa-rate2}
        \gamma=\left(\frac{\sqrt{\kappa^2-1}-1}{\sqrt{\kappa^2-1}+1}\right)^2=1-4\frac{\sqrt{\kappa^2-1}}{\kappa^2+2\sqrt{\kappa^2-1}}.
    \end{equation}
\end{thm}
\begin{pf}
    We note that $\kappa=\sqrt{2}$ implies that $\alpha_1=\alpha_2=1$, otherwise $\alpha_1=\alpha_2\in(1,2)$.
    Convergence to some $x^*\in A\cap B$ follows from convexity.
    If $x^*\not\in \bd A \cap \bd B$, then Theorem~\ref{thm:convex-rate} states that the convergence
    is in finite time, for which the rate holds trivially.
    The remaining case is $x^*\in \bd A \cap \bd B$.
    If $\T_{\bd A}(x^*)=\T_{\bd B}(x^*)$, then $\bd A=\bd B$ in some neighborhood of $x^*$
    and the problem is trivial with convergence in finite time.

    Otherwise, $\T_{\bd A}(x^*)\neq \T_{\bd B}(x^*)$ and consequently the Friedrichs angle satisfies $\cos(\theta_F) > 0$.
    First consider the case where the angle between the sets $A$ and $B$ is obtuse at $x^*$.
    Let $\delta_1$ be such that Lemma~\ref{lem:ball} holds, i.e., $\Proj_A^{\alpha_1} x\in A$ and $\Proj_B^{\alpha_2} x\in B$, for any $x\in\Ball_{\delta_1}(x^*)$.
    Let $c=\langle n_A(x^*), n_B(x^*)\rangle$,
    where $n_A(x^*),n_B(x^*)$ are the outward facing unit normals for the sets $A,B$ at the point $x^*$,
    which by definition of obtuse satisfies $c>0$.
    By smoothness of the boundaries of $A$ and $B$, and continuity of their normals, there is some $\delta_2>0$ such that
    \begin{align}\label{eq:pf:global}
        \langle n_A(x), n_B(y)\rangle>0, \forall x\in\Ball_{\delta_2}(x^*)\cap \bd A, y\in\Ball_{\delta_2}(x^*)\cap \bd B,
    \end{align}
    where $n_A(x)$, $n_B(y)$ are the outward facing unit normals to $A$ and $B$ at $x$ and $y$ respectively.
    Now, by convergence of $x_k$ to $x^*$, there is some $k$ such that $x_k\in\Ball_{\delta}(x^*)$ where $\delta=\min(\delta_1,\delta_2)$.
    Thus by Lemma~\ref{lem:ball} and non-expansiveness of the projectors, we have $\Proj_A^{\alpha_1} x\in A$ and $x_{k+1}=\Proj_B^{\alpha_2}\Proj_A^{\alpha_1} x_k\in B$.
    If $x_{k+1}\in A$, then the problem is solved in finite time, and the result is trivial, otherwise $x_{k+1}\in B\setminus A$.
    There must therefore exist a point $\bar{x}$ on the line between $x_{k+1}\in B\setminus A$ and $\Proj_A^{\alpha_1}x_k\in A$ such that $\bar{x}\in\bd A$,
    moreover it must satisfy $\langle n_A(\bar{x}), x_{k+1} - \Proj_A^{\alpha_1}x_k\rangle > 0$
    since the line is pointing out of the set $A$.
    But by the definition of the projection and $x_{k+1}$, we have
    \begin{align*}
        \frac{x_{k+1} - \Proj_A^{\alpha_1}x_k}{\|x_{k+1} - \Proj_A^{\alpha_1}x_k\|} = -n_B(\tilde{x}),
    \end{align*}
    where $\tilde{x}=\Proj_B\Proj_A^{\alpha_1}x_k\bd B$.
    This leads to $\langle n_A(\bar{x}), n_B(\tilde{x})\rangle <0$.
    Since both $\bar{x}$ and $\tilde{x}$ are in $\Ball_{\delta}(x^*)$ by non-expansiveness,
    this is a contradiction to \eqref{eq:pf:global}, i.e. $x_{x+1}\in B\setminus A$ cannot hold,
    so $x_{x+1}\in A\cap B$ and the convergence is finite and the result holds trivially. 

    The remaining case is when $A$ and $B$ form an acute angle at $x^*$.
    By Lemma~\ref{lem:set-transversality} case \ref{lem:set-transversality:bdAB}, 
    we have $\sreg[A,B](x^*)= \sin(\theta_F/2)\leq 1/\sqrt{2}$, so by definition of $\sreg$ (Definition~\ref{def:sr}),
    it must hold that $\kappa\geq 1/\sreg[A,B](x^*)=1/\sin(\theta_F/2)\geq\sqrt{2}$.
    By Theorem \ref{thm:regularity}, we see that the optimal rate would have been achieved if
    $\kappa=1/\sin(\theta_F/2)$, i.e. $\alpha_1=\alpha_2>\alpha^*$,
    or equivalently that the parameters have been chosen as if $\theta_F$ was smaller.
    But as seen in Remark~\ref{rem:larger-angle},
    this still results in the sub-optimal rate \eqref{eq:kappa-rate2} based on this conservative $\kappa$.
\end{pf}

\begin{rem}
    We note that the adaptive method proposed in \cite{falt-optimal} for estimating $\theta_F$ by the angle between the 
    vectors $v_1=\Proj_B^{\alpha_1}x_k-x_k$ and $v_2=\Proj_A^{\alpha_1}x_k-\Proj_B^{\alpha_2}\Proj_A^{\alpha_1}x_k$,
    works very well in the setting of two convex sets $(A,B)$ with smooth boundaries.
    This can be seen by observing that if $v_1/\|v_1\| = -n_1$ and $v_2/\|v_2\| = n_2$,
    where $n_1,n_2$ are normal vectors with unit length to $A$ and $B$ at the point $x^*$,
    then the angle between them is exactly $\theta_F$ in the acute case.
    And indeed, as long as the algorithm has not already converged, we have $v_1/\|v_1\|\rightarrow -n_1$, $v_2/\|v_2\|\rightarrow n_2$
    as $x_k\rightarrow x^*$, by the definition of the projections and continuity of the normals around $x^*$.
    The estimate will therefore converge to $\theta_F$ as $x_k\rightarrow x^*$.  
\end{rem}

\mnoteh{Something about optimal parameters depending on $x^*$, so choose as a function of $x^*$ if possible. Something about not proven for varying parameters, but convergence if bounded away from $2$.}
\mnoteh{
We now note that when choosing the parameters that are optimal for subspaces, i.e.
\begin{align*}
    \alpha=1,\quad\alpha_1=\alpha_2=\frac{2}{1+\sin(\theta_F)}=:\alpha^*
\end{align*}
where $\theta_F$ if the Friedrichs angle of the tangent spaces at the point $x^*$,
then the convergence for the cases above is given by
\begin{align*}
    \gamma(S_\T(x^*))=\frac{1-\sin(\theta_F)}{1+\sin(\theta_F)}=\alpha^*-1.
\end{align*}
}
\mnoteh{
Optimal alpha vs (c-s)/(c+s):
When optimal reduces by a factor 1e-6 takes more than 100 iterations (up to 100 iterations worse), it is at worst a factor .00239 worse.
When optimal reduces by a factor 1e-6 takes more than 20 iterations (up to 20 iterations worse), it is at worst a best a factor 0.065 worse.
When optimal reduces by a factor 1e-6 takes more than 10 iterations (up to 10 iterations worse), it is at worst a best a factor .0.4 worse.
Conclusion "Optimal" alpha takes no more than 20 iterations or a factor of 6.5\% more iterations whichever is best, than the 2alpha choice.
}

\subsection{Counter example}\label{sec:counter-example}

We present a simple convex example that illustrates that it is not always possible to rely on finite identification of smooth manifolds for the GAP algorithm \ref{eq:GAP},
even in the case of convex polytopes.

\begin{figure}
\centering
\input{cone.tex}
\caption{\label{fig:counter}Illustration of the problem with a cone $C$ and line $D$ from Example \ref{ex:counter}. The iterates $p_0,p_1,p_2,\dots$ are illustrated in red, the normal cone to $C$ with dashed lines, and the rays through $(1,-\gamma)$ and $(-1,-\gamma)$ are shown with blue dotted lines.
As shown in the example, the iterates stay on the dotted lines and alternate between projecting on the two faces of $C$.}
\end{figure}
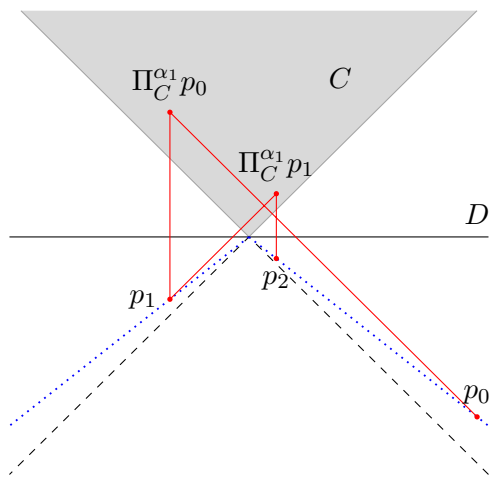

\begin{ex}\label{ex:counter}
Consider the convex feasibility problem ($C,D$) with $C=\left\{ (x,y)\mid y\geq\left|x\right|\right\}$, $D=\left\{ (x,y)\mid y=0\right\}$ as illustrated in Figure \ref{fig:counter},
with parameters $\alpha=1,\alpha_{1}=\alpha_{2}=1.5$ for the GAP algorithm \ref{eq:GAP}. Let
\[
p_{0}=(1,-\gamma)
\]
where $\gamma=\frac{1}{12}\left(1+\sqrt{73}\right)\approx0.795$.
The GAP algorithm will then alternate between projecting onto the half-lines $\{y=x, x>0\}$ and $\{y=-x, x<0\}$.
\end{ex}
\begin{pf}
The first projection point will hit the boundary of the cone $C$ at $\Proj_Cp_0=\frac{1}{2}\left(1-\gamma, 1-\gamma\right)$ which is easily seen by that
$\Proj_Cp_0-p_0=\frac{1}{2}(-1-\gamma, 1+\gamma)\perp\Proj_Cp_0$. The relaxed projection point and the next iterate can then be calculated to
\begin{align*}
\Pi_{C}^{\alpha1}p_{0} & =
\frac{1}{4}\left(1-3\gamma, -3+ \gamma\right)\\
p_{1} & =\Pi_{D}^{\alpha_{2}}\Pi_{C}^{\alpha_{1}}p_{0}=
\frac{1}{8}\left(2-6\gamma, -3+\gamma\right)\\
\end{align*}
We note that $\gamma^2=\frac{1}{6}(\gamma+3)$, and that
$(p_1)_x\gamma=\frac{1}{8}(2-6\gamma)\gamma=\frac{1}{8}(\gamma-3)=(p_1)_y$,
where $p_1=((p_1)_x,(p_1)_y)$.
So $p_{1}$ is simply $p_0$ scaled and flipped around the $y$ axis, i.e., it is of the form $p_{1}=\beta\left(-1,-\gamma\right)$ for some $\beta>0$.
The next projection point is therefore on the boundary of the cone $C$ with $x<0$, and
because of the symmetry around the $y$ axis, the next iterate is
\[
p_{2}=\beta^{2}\left(1,-\gamma\right).
\]

By linearity and induction, it is clear that the algorithm will not identify any
of the smooth surfaces $\{y=x, x>0\}$ or $\{y=-x, x<0\}$ but instead alternate
between them.
\end{pf}

\begin{rem}
    Example~\ref{ex:counter} shows that finite identification of either of the manifolds $\{(x,y) \mid y=x, x>0\}$ and $\{(x,y) \mid y=-x, x<0\}$ does not occur for every initial point.
    However, with some reasonable definition of smallest angle, for example through the subregularity constant $\sreg$, we would have $\theta_F=\pi/4$, and the theory for subspaces would predict a worst case rate $\gamma(S)=0.5$.
    It is notable that the convergence rate $\beta\approx 0.35$ in the example is significantly better.
    It is therefore still an open question whether the smallest angle sets an upper bound on the rate through the eigenvalues in Theorem \ref{thm:gap-all-eigs}
    even for these problems.
\end{rem}

\mnoteh{The optimal parameters achieve a rate of $\beta\approx0.055$ vs the theoretical worst case rate $\approx0.17$.}

\section{Conclusions}
We have shown that the known convergence rates for the GAP algorithm on affine sets extend to local rates on smooth manifolds, and that the optimal parameters and rates hold also in this setting.
These rates are significantly better than previous known rates for similar projection methods.
We have also shown how these results can be applied to generate linear convergence rates for two smooth and solid convex sets, and how they can be connected to linear regularity.

Since finite identification of smooth manifolds cannot generally be assumed, it remains to be shown how these results can be applied to general convex sets.

\appendix
\section{Appendix}

\subsection{Proof of Theorem~\ref{thm:gapP1-1}}\label{app:gapP1-1}
Since $S=T$ with $\alpha=1$,
we begin by showing that all eigenvalues of $T$ in Theorem \ref{thm:gap-all-eigs}
satisfy $|\lambda|\leq\gamma^*$.
For convenience of notation we introduce
\begin{align}
    f(\theta)&\ldef \frac{1}{2}\left(2-\alpha_{1}-\alpha_{2}+\alpha_{1}\alpha_{2}\cos^2(\theta)\right)\\
    g(\theta)&\ldef \sqrt{f(\theta)^{2}-(1-\alpha_{1})(1-\alpha_{2})}
\end{align}
so that $\lambda_{i}^{1,2}$ in \eqref{eq:eig12} can be written $\lambda_{i}^{1,2} = f(\theta_i)\pm g(\theta_i)$.
For $\alpha_1=\alpha_2=\alpha^*=\frac{2}{1+\sin(\theta_F)}$ we get
$f(\theta_{F})=1-\alpha^*+{\alpha^*}^2c_F^2/2=\frac{1-\sin(\theta_F)}{1+\sin(\theta_F)}=\alpha^*-1$ and
$g(\theta_F)=0$.
The eigenvalues corresponding to $\theta_F$ are therefore $\lambda_{F}^{1,2}=\alpha^{*}-1=\frac{1-\sin(\theta_F)}{1+\sin(\theta_F)}$.
We also see that $f(\pi/2)=1-\alpha^{*},\,g(\pi/2)=0$.
Since $f(\theta)$ is linear in $\cos^2(\theta)$, which is decreasing in $\left[\theta_{F},\pi/2\right]$, and $\left|f(\theta_{F})\right|=\left|f(\pi/2)\right|=\alpha^{*}-1$,
it follows that $\left|f(\theta_{i})\right|\leq\alpha^{*}-1$ for all
$\theta_{i}\in\left[\theta_{F},\pi/2\right]$. This means that $f(\theta_{i})^2-(\alpha^{*}-1)^2\leq0$ and the
corresponding $\lambda_{i}^{1,2}$ are complex with magnitudes
\begin{align*}
\left|\lambda_{i}^{1,2}\right| & =\sqrt{f(\theta_{i})^{2}+\left|f(\theta_{i})^{2}-(1-\alpha^{*})^{2}\right|}=\sqrt{(1-\alpha^{*})^{2}}\\
 & =\alpha^{*}-1\quad\forall i:\,\theta_{F}\leq\theta_{i}\leq\pi/2.
\end{align*}
For the remaining eigenvalues we have $|1-\alpha_1|=\alpha^*-1=\gamma^*$,
$|1-\alpha_2|=\alpha^*-1=\gamma^*$,
$|(1-\alpha_1)(1-\alpha_2)|=(\alpha^*-1)^2\leq\gamma^*$.
%
%
Lastly, the eigenvalues in $\lambda=1$ correspond to the angles $\theta_{i}=0$ and are
semisimple since the matrix in~\eqref{eq:T1matrix} is diagonal for $\theta_{i}=0$.
We therefore conclude, using
Fact~\ref{fct:limitexists} and \ref{fact:Convergent-at-rate},
that $\alpha_{1}=\alpha_{2}=\alpha^{*}$ results in that the GAP operator $S=T$ in
\eqref{eq:GAP2} is linearly convergent with any rate $\mu\in\left(\gamma^{*},1\right)$
where $\gamma^{*}=\alpha^{*}-1=\frac{1-\sin(\theta_F)}{1+\sin(\theta_F)}$
is a subdominant eigenvalue.

\subsection{Lemmas}

\label{app:lemma}
\begin{lem}[Infinite Sub-sequence]\label{lem:sub-sequence}%
    Given any infinite sequence of increasing positive integers
    $(r_{j})_{{j\in\naturalN}}$, for any integer $n>0$ there exists an infinite sub-sequence $(r_{j_{k}})_{k\in\naturalN}$ where
    \[
    r_{j_{k}}=a+nb_{k},
    \]
    for some $a\in\naturalN$ and some increasing sequence $(b_{k})_{k\in\naturalN}$.
\end{lem}
\begin{pf}
    Fix $n$ and consider the finite collection of sets $S_i=\{v\in\naturalN\mid v = i+nb, b\in\naturalN\}$,
    $i=0,\ldots,n-1$.
    We have $\cup_{i=0,\ldots,n-1}S_i=\naturalN$,
    so $\cup_{i=0,\ldots,n-1}(S_i\cap \{r_j\}_j)=\{r_j\}_{j\in\naturalN}$
    and thus one of the sets $(S_i\cap \{r_j\}_{j\in\naturalN})$ must be infinite.
    Let $a$ be the index so that $(S_a\cap \{r_j\}_{j\in\naturalN})$ is infinite.
    This is clearly a subset of $\{r_j\}_{j\in\naturalN}$ and by the definition of $S_a$
    each element is of the form $a+nb_{k}$ with $b_{k}\in\naturalN$ and the proof is complete.
\end{pf}

\begin{lem}\label{lem:trdet}
The matrix
\begin{align}\label{eq:lemma:M}
    M\ldef (2-\alpha^{*})I+\frac{\alpha^{*}}{\alpha_{1}}(T_{1}^{F}-I),
\end{align}
where $T_{1}^{F}$ is the matrix defined in~\eqref{eq:T1matrix} corresponding to the angle $\theta_{F}$
has trace and determinant:
\begin{align*}
\text{tr}M & =  \frac{2}{(1+s)\alpha_{1}}\left(-\alpha_{1}-\alpha_{2}+\alpha_{2}\alpha_{1}c^{2}+2\alpha_{1}s\right)\\
\det M & =  \frac{4s(1-s)}{\alpha_{1}(1+s)^{2}}\left(-\alpha_{1}-\alpha_{2}+\alpha_{1}\alpha_{2}(1+s)\right),
\end{align*}
where $s\ldef \sin(\theta_F),\,c\ldef \cos(\theta_F)$.
\end{lem}

\begin{pf}
The matrix $M$ can be written
\begin{align*}
M &= (2-\alpha^{*})I+\frac{\alpha^{*}}{\alpha_{1}}\left(\begin{pmatrix}1-\alpha_{1}s^{2} & \alpha_{1}cs\\
\alpha_{1}(1-\alpha_{2})cs & (1-\alpha_{2})(1-\alpha_{1}c^{2})
\end{pmatrix}-I\right)\\
 &= \begin{pmatrix}2-\alpha^{*}-\alpha^{*}s^{2} & \alpha^{*}cs\\
\alpha^{*}(1-\alpha_{2})cs & 2-\alpha^{*}+\frac{\alpha^{*}}{\alpha_{1}}\left((1-\alpha_{2})(1-\alpha_{1}c^{2})-1\right)
\end{pmatrix}\\
 &= \begin{pmatrix}2-\alpha^{*}(1+s^{2}) & \alpha^{*}cs\\
\alpha^{*}(1-\alpha_{2})cs & 2-\alpha^{*}+\frac{\alpha^{*}}{\alpha_{1}}\left(\alpha_{1}\alpha_{2}c^{2}-\alpha_{2}-\alpha_{1}c^{2}\right)
\end{pmatrix}.
\end{align*}
Using that $\alpha^{*}=\frac{2}{1+s}$, we can rewrite the diagonal elements
\[
2-\alpha^{*}(1+s^{2})=\alpha^{*}\left(1+s-(1+s^{2})\right)=\alpha^{*}s(1-s)
\]
and
\begin{align*}
2-\alpha^{*}+\frac{\alpha^{*}}{\alpha_{1}}\left(\alpha_{1}\alpha_{2}c^{2}-\alpha_{2}-\alpha_{1}c^{2}\right)
&= \alpha^{*}(1+s)-\alpha^{*}+\alpha^{*}\left(c^{2}(\alpha_{2}-1)-\frac{\alpha_{2}}{\alpha_{1}}\right)\\
 & = \alpha^{*}\left(s+c^{2}(\alpha_{2}-1)-\frac{\alpha_{2}}{\alpha_{1}}\right).
\end{align*}
We can extract the factor $\alpha^{*}cs$ from the matrix and get
\[
M=\alpha^{*}cs\begin{pmatrix}\frac{1-s}{c} & 1\\
1-\alpha_{2} & \frac{s+c^{2}(\alpha_{2}-1)-\frac{\alpha_{2}}{\alpha_{1}}}{cs}
\end{pmatrix}.
\]
The trace is therefore given by
\begin{align*}
\text{tr}M & = \alpha^{*}cs\left(\frac{1-s}{c}+\frac{s+c^{2}(\alpha_{2}-1)-\frac{\alpha_{2}}{\alpha_{1}}}{cs}\right)\\
 & = \alpha^{*}\left(2s-s^{2}+c^{2}\alpha_{2}-c^{2}-\frac{\alpha_{2}}{\alpha_{1}}\right)\\
 & = \frac{\alpha^{*}}{\alpha_{1}}\left(-\alpha_{1}-\alpha_{2}+\alpha_{2}\alpha_{1}c^{2}+2\alpha_{1}s\right)\\
 & = \frac{2}{(1+s)\alpha_{1}}\left(-\alpha_{1}-\alpha_{2}+\alpha_{2}\alpha_{1}c^{2}+2\alpha_{1}s\right)
\end{align*}
and the determinant is given by
\begin{align*}
\text{det}M & = \left(\alpha^{*}cs\right)^{2}\left(\frac{\left(1-s\right)\left(s+c^{2}(\alpha_{2}-1)-\frac{\alpha_{2}}{\alpha_{1}}\right)}{c^{2}s}-\frac{\left(1-\alpha_{2}\right)c^{2}s}{c^{2}s}\right)\\
& = \alpha^{*2}s\biggl(\left(s+c^{2}(\alpha_{2}-1)-\frac{\alpha_{2}}{\alpha_{1}}-s^{2}-c^{2}s(\alpha_{2}-1)+s\frac{\alpha_{2}}{\alpha_{1}}\right)-\left(1-\alpha_{2}\right)c^{2}s\biggr)\\
& = \alpha^{*2}s\left(s+c^{2}(\alpha_{2}-1)-\frac{\alpha_{2}}{\alpha_{1}}-s^{2}+s\frac{\alpha_{2}}{\alpha_{1}}\right)\\
& = \alpha^{*2}s\left(s-1+\alpha_{2}c^{2}+\frac{\alpha_{2}}{\alpha_{1}}(s-1)\right)\\
& = \alpha^{*2}s(1-s)\left(-1+\alpha_{2}(1+s)-\frac{\alpha_{2}}{\alpha_{1}}\right)\\
& = \frac{\alpha^{*2}s(1-s)}{\alpha_{1}}\left(-\alpha_{1}-\alpha_{2}+\alpha_{1}\alpha_{2}(1+s)\right)\\
& = \frac{4s(1-s)}{\alpha_{1}(1+s)^{2}}\left(-\alpha_{1}-\alpha_{2}+\alpha_{1}\alpha_{2}(1+s)\right).
\end{align*}
\end{pf}

\begin{lem}\label{lem:realpart}
Under the assumptions $\alpha=\frac{\alpha^{*}}{\alpha_{1}}$, $\cmt{0\neq}\alpha_{1}\geq\alpha_{2}\cmt{\neq0}>0$
and $\theta_F\in(0,\pi/2)$, the matrix $M$ \eqref{eq:lemma:M} in Lemma~\ref{lem:trdet} satisfies
\[
\left(\alpha_{1}\neq\alpha^{*}\text{ or }
\alpha_{2}\neq\alpha^{*}\right)
\Rightarrow\max\re\spectrum{M}>0,
\]
where $\Lambda(M)$ is the set of eigenvalues of $M$.
\end{lem}
\begin{pf}
We prove the equivalent claim
\[
\max\re\spectrum{M}\leq0\Rightarrow\alpha_{1}=\alpha_{2}=\alpha^{*}.
\]
We have $\max\re\spectrum{M}\leq0$ if and only if both eigenvalues of $M$
have negative or zero real part, which is equivalent to
\[
\lambda_{1}+\lambda_{2}\leq0\quad\text{and}\quad\lambda_{1}\lambda_{2}\geq0.
\]
This is equivalent to
\[
\text{tr}M\leq0\quad\text{and}\quad\text{det}M\geq0.
\]
Using Lemma~\ref{lem:trdet}, this can be written
\begin{align*}
\begin{cases}
\frac{2}{(1+s)\alpha_{1}}\left(-\alpha_{1}-\alpha_{2}+\alpha_{2}\alpha_{1}c^{2}+2\alpha_{1}s\right) & \leq 0\\
\frac{4s(1-s)}{\alpha_{1}(1+s)^{2}}\left(-\alpha_{1}-\alpha_{2}+\alpha_{1}\alpha_{2}(1+s)\right) & \geq 0
\end{cases},
\end{align*}
where $s\ldef\sin(\theta_F)$ and $c\ldef\cos(\theta_F)$.
Since $\alpha_1>0$, $s\in(0,1)$, this is equivalent to
\begin{subnumcases}{\label{eq:maxrem0}}
   \cmt{\sign{\alpha_1}(}\alpha_{1}+\alpha_{2}-\alpha_{2}\alpha_{1}c^{2}-2\alpha_{1}s\cmt{)} & $\geq 0$\label{eq:maxrem}\\
   \cmt{\sign{\alpha_1}(}-\alpha_{1}-\alpha_{2}+\alpha_{1}\alpha_{2}(1+s)\cmt{)} & $\geq 0$.\label{eq:maxrem2}
\end{subnumcases}
This implies that the sum is positive, i.e.
\begin{align*}
\cmt{\sign{\alpha_1}}\big(\alpha_{1}+\alpha_{2}-\alpha_{2}\alpha_{1}c^{2}&-2\alpha_{1}s\big)+
\cmt{\sign{\alpha_1}}\left(-\alpha_{1}-\alpha_{2}+\alpha_{1}\alpha_{2}(1+s)\right)\\
& = \cmt{\sign{\alpha_1}}(\alpha_{2}\alpha_{1}s^{2}-2\alpha_{1}s+\alpha_{1}\alpha_{2}s)\\
& = \cmt{|}\alpha_{1}\cmt{|}s\left(\alpha_{2}s-2+\alpha_{2}\right) \geq 0
\end{align*}
which, since $\alpha_2,s>0$, is equivalent to $\alpha_{2}(1+s) \geq 2$, and thus
\begin{align*}
\alpha_{2} \geq \frac{2}{1+s}=\alpha^{*}.
\end{align*}
But since $\alpha_{2}\geq\alpha^{*}$\cmt{, by assumption $\alpha_1\geq\alpha_2$ and thus $\sign{\alpha_1}=1$, and},~\eqref{eq:maxrem} implies
\[
\alpha_{1}+\alpha_{2}-\alpha^{*}\alpha_{1}c^{2}-2\alpha_{1}s\geq0
\]
which is equivalent to
\begin{align*}
\alpha_{1}+\alpha_{2}-\alpha^{*}\alpha_{1}c^{2}-2\alpha_{1}s&=\alpha_{1}+\alpha_{2}-2\alpha_{1}(1-s)-2\alpha_{1}s\\
&=\alpha_{1}+\alpha_{2}-2\alpha_{1}= \alpha_{2}-\alpha_{1} \geq 0
\end{align*}
i.e., $\alpha_{2}\geq\alpha_{1}.$

But by the assumption that $\alpha_{1}\geq\alpha_{2}$ we know that \eqref{eq:maxrem0}
implies $\alpha_{1}=\alpha_{2}\geq\alpha^{*}$. Equation~\eqref{eq:maxrem}
yields
\begin{alignat*}{2}
&& \alpha_{1}+\alpha_{2}-\alpha_{2}\alpha_{1}c^{2}-2\alpha_{1}s &\geq0\\
&\Rightarrow \quad&
2\alpha_{1}-\alpha_{1}^{2}c^{2}-2\alpha_{1}s&\geq0\\
 &\Leftrightarrow &
2-\alpha_{1}c^{2}-2s&\geq0\\
&\Leftrightarrow & 2\frac{(1-s)}{c^{2}}&\geq\alpha_{1}\\
&\Leftrightarrow & \alpha^*=\frac{2}{(1+s)}&\geq\alpha_{1},
\end{alignat*}
where the implication is from $\alpha_{1}=\alpha_{2}$.
We have shown that $\alpha^{*}\geq\alpha_{1}=\alpha_{2}\geq\alpha^{*}$
i.e., $\alpha^{*}=\alpha_{1}=\alpha_{2}\geq\alpha^{*}$. This completes the proof.
\end{pf}

\subsection{Proof of Theorem~\ref{thm:gapiff}}\label{app:gapiff}
The first direction, that both $S_1$ and $S_2$ are convergent with any rate $\mu\in(\gamma^*,1)$ for the parameters in \eqref{eq:optpar} holds by Theorem~\ref{thm:gapP1-1}.
We now prove that if $S_1$ and $S_2$ converge with the rate $\mu$ for all $\mu\in(\gamma^*,1)$ then the parameters must be those in \eqref{eq:optpar}.
By Fact \ref{fct:limitexists}, if both operators converge with any rate $\mu\in(\gamma^*,1)$ then it must be that
$\gamma(S_1)\leq \gamma^*$ and $\gamma(S_2)\leq \gamma^*$.
By Definition \ref{def:subdominant}, this means that all eigenvalues $\lambda$ to both $S_1$ and $S_2$
have $|\lambda|\leq \gamma^*$, unless $\lambda=1$.
With $S_i=(1-\alpha)I + \alpha T_i$, we see from Theorem \ref{thm:gap-all-eigs},
that $T_1$ has an eigenvalue in $1-\alpha_2$, $T_2$ in $1-\alpha_1$,
and both $T_1$ and $T_2$ have eigenvalues in $\lambda_i^{1,2}$ corresponding to the angle $\theta_F$.
We therefore need that $|1+\alpha\left(\lambda-1\right)|\leq\gamma^*$ for each of the eigenvalues $\lambda$.
We start by defining $\hat{\alpha}=\alpha^*/\alpha_1$, where $\alpha^*=2/(1+\sin(\theta_F))$,
and observe that $\alpha^*-1=\gamma^*$.

Assume that $\alpha_1\geq\alpha_2$ and $\alpha=\hat{\alpha}$.
For the eigenvalue $\lambda=1-\alpha_{1}$, we get
\begin{align}\label{eq:eigin1mina}
1+\hat{\alpha}(\lambda-1)=1+\frac{\alpha^{*}}{\alpha_{1}}(1-\alpha_{1}-1)=1-\alpha^{*}.
\end{align}
Consider the eigenvalues of $I+\hat{\alpha}(T_F-I)$ where $T_F$
is the matrix~\eqref{eq:T1matrix} corresponding to the angle $\theta_{F}$,
i.e., the eigenvalues $\lambda_i^{1,2}$.
We have
\begin{equation}
\max\re\spectrum{I+\hat{\alpha}(T_F-I)}>\alpha^{*}-1\label{eq:maxrealpart}
\end{equation}
if and only if
\begin{equation}
\max\re\spectrum{(2-\alpha^{*})I+\hat{\alpha}(T_F-I)}>0.\label{eq:eigpositive}
\end{equation}

By Lemma~\ref{lem:realpart} we know that~\eqref{eq:eigpositive} is true when $\alpha=\hat{\alpha}$, unless $\alpha_1=\alpha_2=\alpha^*$.
We therefore know that for $\alpha=\hat\alpha$, unless the optimal parameters are selected,
there will always be\cmt{~one case (C1), with} one eigenvalue of $S_2$ in
$1-\alpha^*$ and one\cmt{~eigenvalue}, corresponding to $\theta_F$, with real part greater than $\alpha^*-1$.
We now consider the two cases $\alpha>\hat{\alpha}$ and $\alpha<\hat{\alpha}$. First note that $\alpha$ acts as a scaling of the eigenvalues relative to the point $1$, i.e., $(1-\alpha)+\alpha\lambda=1+\alpha(\lambda-1)$.
It is therefore clear that $\alpha>\hat\alpha$ will result in one eigenvalue with real part less than $1-\alpha^*=-\gamma^*$, and thus $\gamma(S_1)>\gamma^*$ and $\gamma(S_2)>\gamma^*$.

Similarly, any $\alpha<\hat\alpha$ will result in one eigenvalue ($\lambda_F^1$) with real part greater than $\alpha^*-1=\gamma^*$. If this eigenvalue is not in $1$, i.e., unless $1+\alpha(\lambda_F^1-1)=1$, we know that $\gamma(S)>\gamma^*$ also in this case.
Since $\alpha\neq0$ we have $1+\alpha(\lambda_F^1-1)=1$ if and only if $\lambda_F^1=1$.
But $\lambda_F^1=1$ only if $\det(T_F-I)=0$,
where $T_F$ is the block corresponding to $\theta_F$ in \eqref{eq:T1matrix}.
Since $\alpha_1,\alpha_2\neq0$ and $\theta_F>0$ we get
\begin{align*}
    \det(T_F-I)=
-\alpha_1s_F^2(\alpha_1c_F^2-\alpha_2+\alpha_1\alpha_2c_F^2)-\alpha_1^2(1-\alpha_2)c_F^2s_F^2=
\alpha_1\alpha_2s_F^2\neq0
\end{align*}
and thus $\lambda_F^1\neq1$.

We conclude that when $\alpha_1\geq\alpha_2$, then $\gamma(S_2)>\alpha^*-1=\gamma^*$
for all parameters that are not $\alpha=1, \alpha_1=\alpha_2=\alpha^*$.

The proof is only dependent on the eigenvalue $1-\alpha_1$, corresponding to $S_2$,
and the eigenvalue $\lambda_F^{1,2}$ corresponding to $\theta_F$.
From symmetry of $\alpha_{1},\alpha_{2}$ in $\lambda_F^{1,2}$ we
see that the same argument holds if we instead assume $\alpha_{2}\geq\alpha_{1}$, let $\hat\alpha=\alpha^*/\alpha_2$,
and consider the eigenvalues $1-\alpha_2$ from $S_1$ and $\lambda_F^{1,2}$.
Therefore, when $\alpha_2\geq\alpha_1$, we get $\gamma(S_1)>\alpha^*-1=\gamma^*$
for all parameters that are not $\alpha=1, \alpha_1=\alpha_2=\alpha^*$.
To conclude, unless $\alpha=1,\alpha_1=\alpha_2=\alpha^*$,
we have either $\gamma(S_1)>\gamma^*$ or $\gamma(S_2)>\gamma^*$,
which contradicts that they both converge linearly with any rate $\mu\in(\gamma^*,1)$.

\bibliography{references}
\bibliographystyle{jnsao}

\end{document}

%% file: cone.tex
\begin{tikzpicture}[scale=3]

\draw[fill=gray,opacity=0.3]
  (1,1) -- (0,0) -- (-1,1);
\node at (0.4,0.7) {$C$};
\draw[dashed] (1.05,-1.05) -- (0,0) -- (-1.05,-1.05);
\coordinate (g1) at (1, -0.7953336454431275);
\coordinate (g2) at (-1, -0.7953336454431275);

\draw (-1.05, 0) -- (1.05,0);
\node at (1.0,0.1) {$D$};

\draw[dotted,blue,line width=0.7] (0,0) -- ($1.05*(g1)$);
\draw[dotted,blue,line width=0.7] (0,0) -- ($1.05*(g2)$);

\draw[red]
    (g1) node[fill=red,circle,outer sep=0.0, inner sep=0.7,label={[color=black]above:$p_0$}] {} --
    (-0.3465,0.5511) node[fill=red,circle,outer sep=0.0, inner sep=0.7, label={[color=black]above:$\Pi_C^{\alpha_1}p_0$}] {} --
    (-0.3465, -0.2755) node[fill=red,circle,outer sep=0.0, inner sep=0.7,label={[color=black]left:$p_1$}] {} --
    (0.1200623311095328, 0.19095627900278075) node[fill=red,circle,outer sep=0.0, inner sep=0.7,label={[color=black]above:$\Pi_C^{\alpha_1}p_1$}] {} --
    ($0.12*(g1)$) node[fill=red,circle,outer sep=0.0, inner sep=0.7,label={[color=black]below:$p_2$}] {};

\end{tikzpicture}